\documentclass[graybox]{svmult}

\usepackage{mathptmx}       % selects Times Roman as basic font
\usepackage{helvet}         % selects Helvetica as sans-serif font
\usepackage{courier}        % selects Courier as typewriter font
\usepackage{type1cm}        % activate if the above 3 fonts are
\usepackage{makeidx}         % allows index generation
\usepackage{graphicx}        % standard LaTeX graphics tool
\usepackage{multicol}        % used for the two-column index
\usepackage[bottom]{footmisc}% places footnotes at page bottom
\usepackage{caption}
\usepackage{subcaption}

\spdefaulttheorem{project}{Research Project}{\bf}{}
\newcommand{\resproject}[1]{\begin{svgraybox} \begin{project} #1 \end{project} \end{svgraybox}}

\spdefaulttheorem{cproblem}{Challenge Problem}{\bf}{}
\newcommand{\chproblem}[1]{\begin{cproblem} #1 \end{cproblem}}

\spnewtheorem*{prerequisites}{Suggested prerequisites}{\bf}{\small\em}

\makeindex             % used for the subject index
\begin{document}

\title*{Using Integer Programming to Solve Games, Puzzles, and Ciphers}
% Use \titlerunning{Short Title} for an abbreviated version of
% your contribution title if the original one is too long
\author{Elizabeth L. Bouzarth, John M. Harris, Kevin R. Hutson, and Christian R. Millichap}
\authorrunning{E.L. Bouzarth, J.M. Harris, K.R. Hutson, and C. Millichap}
% Use \authorrunning{Short Title} for an abbreviated version of
% your contribution title if the original one is too long
\institute{Elizabeth L. Bouzarth \at Furman University, 3300 Poinsett Hwy., Greenville, SC 29613, \email{liz.bouzarth@furman.edu}
\and John M. Harris \at Furman University, 3300 Poinsett Hwy., Greenville, SC 29613, \email{john.harris@furman.edu}
\and Kevin R. Hutson \at Furman University, 3300 Poinsett Hwy., Greenville, SC 29613, \email{kevin.hutson@furman.edu}
\and Christian R. Millichap \at Furman University, 3300 Poinsett Hwy., Greenville, SC 29613, \email{christian.millichap@furman.edu}}
%
% Use the package "url.sty" to avoid
% problems with special characters
% used in your e-mail or web address
%
\maketitle

\abstract*{In this paper, we introduce three different classes of undergraduate research projects that implement model building and integer programming. These research projects focus on determining and analyzing solutions to the puzzle \emph{The Genius Square}, optimizing allocation of trains to maximize points in the game \emph{Ticket to Ride}, and (code)breaking monoalphabetic substitution ciphers. Initial models and analysis that came from undergraduate research projects that the authors supervised are shared along with  a variety of open research questions.}

\abstract{In this paper, we introduce three different classes of undergraduate research projects that implement model building and integer programming. These research projects focus on determining and analyzing solutions to the game \emph{The Genius Square}, optimizing allocation of trains to maximize points in the game \emph{Ticket to Ride}, and (code)breaking monoalphabetic substitution ciphers. Initial models and analyses for these scenarios that came from previous undergraduate research projects are shared along with  a variety of open research questions.}

\begin{prerequisites}
Linear Algebra, Introduction to Programming (Python, R, etc.)

\end{prerequisites}

\section{Introduction:  Integer Programming}
\label{sec:1}
%Use the template \emph{author.tex} together with the Springer document class SVMono (monograph-type books) or SVMult (edited books) to style the various elements of your chapter content in the Springer layout.  
%Instead of simply listing headings of different levels we recommend to
%let every heading be followed by at least a short passage of text.
%Further on please use the \LaTeX\ automatism for all your
%cross-references and citations. And please note that the first line of
%text that follows a heading is not indented, whereas the first lines of
%all subsequent paragraphs are.

A mathematical model can be helpful when one wants to represent the structure, features, and function of some object or process.  This representation can often be used to study internal relationships of systems and to make decisions based on these relationships.  Integer programming (IP) is a class of constrained optimization problems where modelers set up integer-valued decision variables, say $x_1, x_2 ,\ldots, x_n$, to encode possible solutions to questions such as how many units of an item should be manufactured/shipped/ordered or where a new fire station, retail complex, or post office should be located to best serve a community.  In addition to setting up \emph{decision variables}, the modeler sets up an \emph{objective function} and \emph{constraints}, which involve these decision variables and represent, respectively, the goal and limiting scenario factors that are important in the decision.  In many cases, the constraints and objective function are required to be linear expressions of the decision variables.  Models like this
are usually expressed in the following manner:

\noindent
\begin{center}
\begin{tabular} {lclcl}
maximize & &$\displaystyle z= \sum_{j} c_j x_j$& & \\
subject to & &$\displaystyle \sum_j a_{ij}x_j \leq b_i$ && $i=1, 2, \ldots, m$ \\
 && $x_j \geq 0$ (integer) & &$j = 1, 2, \ldots, n$. \\
\end{tabular}
\end{center} 
Here, $z$ represents the objective function, and the inequalities that follow ``subject to" are the constraints. The goal of such a model is to find the values of the decision variables ($x_i$) that maximize the value of $z$, the objective function, while satisfying all of the constraints.  Models such as this are often
used to make decisions in a variety of industries.  Examples include supply chain management, production/inventory control, logistics and transportation, sports league scheduling \cite{ams}, and donor/patient transplant matching \cite{Gentry}.  See \cite{HillLieb} for more examples.  

The area of integer programming is also a nice playground in which undergraduates can explore modeling decisions.  While some extensive problems that arise in industry may be more appropriate for graduate work, we have found plenty of applications that are interesting to undergraduate students.  Many popular strategy games and puzzles require players to make decisions, and there is an industry built around helping players make better decisions in these games (e.g. \cite{Bettles, Reiber,Slade,WitterLyford}).  In this paper, we show how we have used integer programming with undergraduate students to model and produce solutions to several popular games, puzzles, and ciphers.  In each case the students involved had no prior experience with the  concepts of integer programming.  They were able to pick up on the modeling concepts, learn how to program a solver, and analyze solutions within an eight-week summer research program.  Each of the puzzles we describe can be modeled with binary variables, and they make use of logic-type constraints to produce solutions.  The focus of these student projects has been to build and evaluate IP models by determining the appropriate variables and constraints while testing different objective functions. Students learn how to use an industrial solver (e.g., CPLEX, Gurobi, GLPK) which implements the best IP algorithms to produce solutions to the models. For information on getting started using optimization solvers, see \cite{But, Maarefdoost, Spokutta}.  For an undergraduate, we believe there is value and creativity in the model-creation process and in the discipline it takes to write and debug code.  Focusing on model-building also minimizes background requirements, and this often allows less experienced students to be involved in research.  Students are often inspired by this process to take an Operations
Research course to learn more about model building and about algorithms used to solve the models.

Before introducing  three applications we explored with undergraduate researchers in Sections \ref{sec:3}-\ref{sec:5}, we present some more traditional examples of building IP models for the reader in Section \ref{sec:2}.

\section{Integer Programming Examples}
\label{sec:2}

Suppose we manufacture two products, $A$ and $B$, which each require two resources, $C$ and $D$.  For every unit we sell of $A$, we make \$5 in profit, and for every unit we sell of $B$ we make \$7 in profit.  Further, we need three units of $C$ and five units of $D$ to make a unit of product $A$, and we need two units of $C$ and eight units of $D$ to make a unit of product $B$.  Suppose also that we only have 100 units of $C$ and 50 units of $D$ in inventory.  We wish to find the right mix of $A$ and $B$ to produce that will maximize our profits.  

To model this, we set up decision variables, $x_A$ and $x_B$, to encode how many units of $A$ and $B$ to make, respectively.  These variables must have integer values since producing fractional amounts of $A$ or $B$ is not appropriate.   The integer programming model is then as follows:

\noindent
\begin{center}
\begin{tabular} {lclcl}
maximize && $ z= 5x_A + 7x_B$ && \\
subject to && $3x_A + 2x_B \leq 100$ &&(inventory $C$ constraint) \\
 && $5x_A + 8x_B \leq 50$  &&(inventory $D$ constraint)\\
 && $x_A, x_B \geq 0$ &&(integer). \\
\end{tabular}
\end{center}
Here, $z$ represents our profit function, which combines the profit made from selling the $A$ units with the profit made from selling the $B$ units.  The two inventory constraints  say that the number of units of $C$ and $D$ that are used in manufacturing $A$ and $B$ cannot exceed the available units of inventory for $C$ and $D$.  Finally, the last constraint ensures that production cannot be negative and must be integer-valued.  
While the context of this example is a manufacturing process, it should be noted that the same idea applies in many scenarios where scarce resources need to be
allocated in order to meet some objective.

Integer programming is also used to model scenarios that require yes/no decisions.  For example, a workforce manager may need to decide whether to assign a certain employee to a certain job, or a university housing director may need to decide whether to assign a pair of roommates to a particular room.  In either case, the yes/no decision can be modeled with binary variables, $x_{ij}$.  For example, we might say that $x_{ij}=1$ if we decide,  `Yes, assign roommate pair $i$ to room $j$', and $x_{ij}=0$ if we decide `No.'    These binary (Boolean) variables then can be used in constraints to make sure that certain rules are followed (for instance, at most one roommate pair can be assigned to each room).

In many  introduction to proof and discrete mathematics courses, students learn rules of Boolean logic, combining two or more statements into compound statements and determining in what scenarios those statements are true or false.  These scenarios can be modeled with constraints on binary variables in IP.  For a more extensive treatment of logical statements modeled as constraints see Chapter 12 of \cite{HillLieb} or \cite{ForWad}.

Suppose we wanted to enforce the conjunction of two variables.  Conjunctions ($A \cap B$) are true only if both statements $A$ and $B$ are true.  In IP, we establish a constraint to ensure both binary decision variables $x_A$ and $x_B$ are assigned the value 1 with the following constraint:
$$x_A + x_B = 2.$$
In modeling the disjunction ($A \cup B$), we seek a constraint that ensures the statement is true when either $A$ or $B$ is true.  The corresponding linear constraint is 
$$x_A + x_B \geq 1.$$
 Finally, we can model implication statements  like `if $A$, then $B$'.  Here, if $A$ is true then $B$ must be true, and if $A$ is false, then $B$ can be either true or false.  The corresponding linear constraint that ensures this is
$$x_A \leq x_B.$$
Each of these can be extended to more variables.  For example, the constraints 
\[
\begin{array}{l}
x_1 + x_2 + \cdots + x_n \geq 1, \\
x_1 + x_2 + \cdots + x_n = 1,\mbox{ and} \\
x_1 + x_2 + \cdots + x_n \leq 1,
\end{array}
\]
\noindent
model that at least, exactly, or at most one of the binary variables $x_1, \ldots, x_n$ can be chosen to have the value of 1, respectively.  Further, the constraints
\[
\begin{array}{l}
x_1 \leq x_2, \\
x_1 \leq x_3, \mbox{ and} \\
x_1 \geq x_2 + x_3 -1
\end{array}
\]
\noindent
ensure that $x_1$ takes on the value 1 if and only if both $x_2$ and $x_3$ take on the value 1.

\begin{exercise}
If $x_A \leq x_B$ models the case where $x_B=1$ if $x_A=1$, what constraint might model $x_B=1$ if and only if $x_A=1$?
\end{exercise}

%In the remaining sections of this chapter, we explore three types of problems we have used IP and associated logic to solve with undergraduate research students. 
In the remaining sections of this chapter, we explore three types of problems we have worked on with undergraduate students that implement IP and associated logic. The setup of these problems are independent of each other, so the reader can choose to read Sections \ref{sec:3}-\ref{sec:5} in any order or depth they wish without compromising their ability to continue working through the chapter.

%%%%%%%%%%%%%%%%%%%%%%%%%%%%%%%%%%%%%%%%%%%%%%%%%%%%%%%%%%%%%%%%%%%%%%%%%%%%%%%%%%%%%%%%%%%%%%%%%%%
  
\section{Genius Square}
\label{sec:3}
% Always give a unique label
% and use \ref{<label>} for cross-references
% and \cite{<label>} for bibliographic references
% use \sectionmark{}
% to alter or adjust the section heading in the running head

%Genius Square is a one-player game where players create a puzzle to solve.  The puzzle consists of a $6 \times 6$ grid game board, seven six-sided game dice, seven cylinder ``blockers", and nine polyominoe pieces.  To create a puzzle, the player rolls the seven dice to reveal locations to place the seven cylinders on the game board.  Each die has a number of game grid locations on its faces.  Figure \ref{Dice} shows the grid board locations of each die.  For example, as seen in Figure \ref{Dice}, one die, shown in red, has the locations $A2$, $A3$, $B1$, $B2$, $B3$, and $C2$ on its faces, where the grid rows are labeled $A$-$F$ from top to bottom and the columns are labeled $1-6$ left to right.  Another die, shown in purple in Figure \ref{Dice}, has only the locations $A5$, $B6$, $E1$ and $F2$ over its six faces.   Figure \ref{GSEx} shows an example of the placement of the blockers based on a roll of the dice.  \christian{Once we get to the decision variables section, we switch to $(i,j)$ coordinates where $1 \leq i,j \leq 6$. Unless using the (letter, number) coordinates serve a purpose, I would just use the $(i,j)$ coordinates from the start and avoid a transition.} \kevin{We use the alphabet-number notation first because that is what is on the game dice.  It is easier to model summations with numbers which insprires the switch.  I'm okay being consistent though.}

\begin{figure}[h!]
	\begin{center}
		\includegraphics[width=0.65 \textwidth]{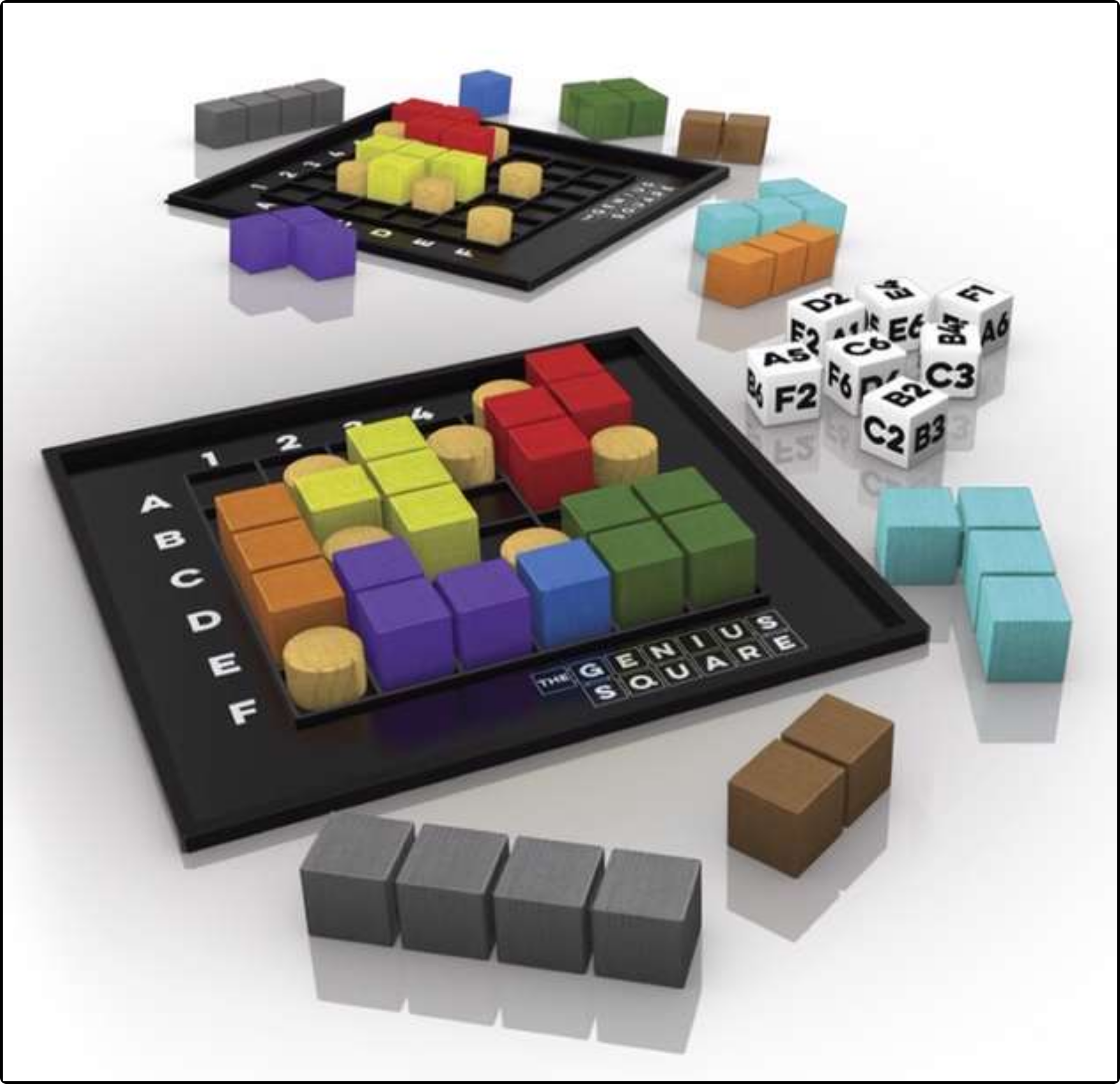}
	\end{center}
	\caption{\textit{The Genius Square}, image courtesy of \cite{hpc}.}
	\label{GS}
\end{figure}

\begin{figure}[h!]
	\begin{center}
		\includegraphics[width=0.35 \textwidth]{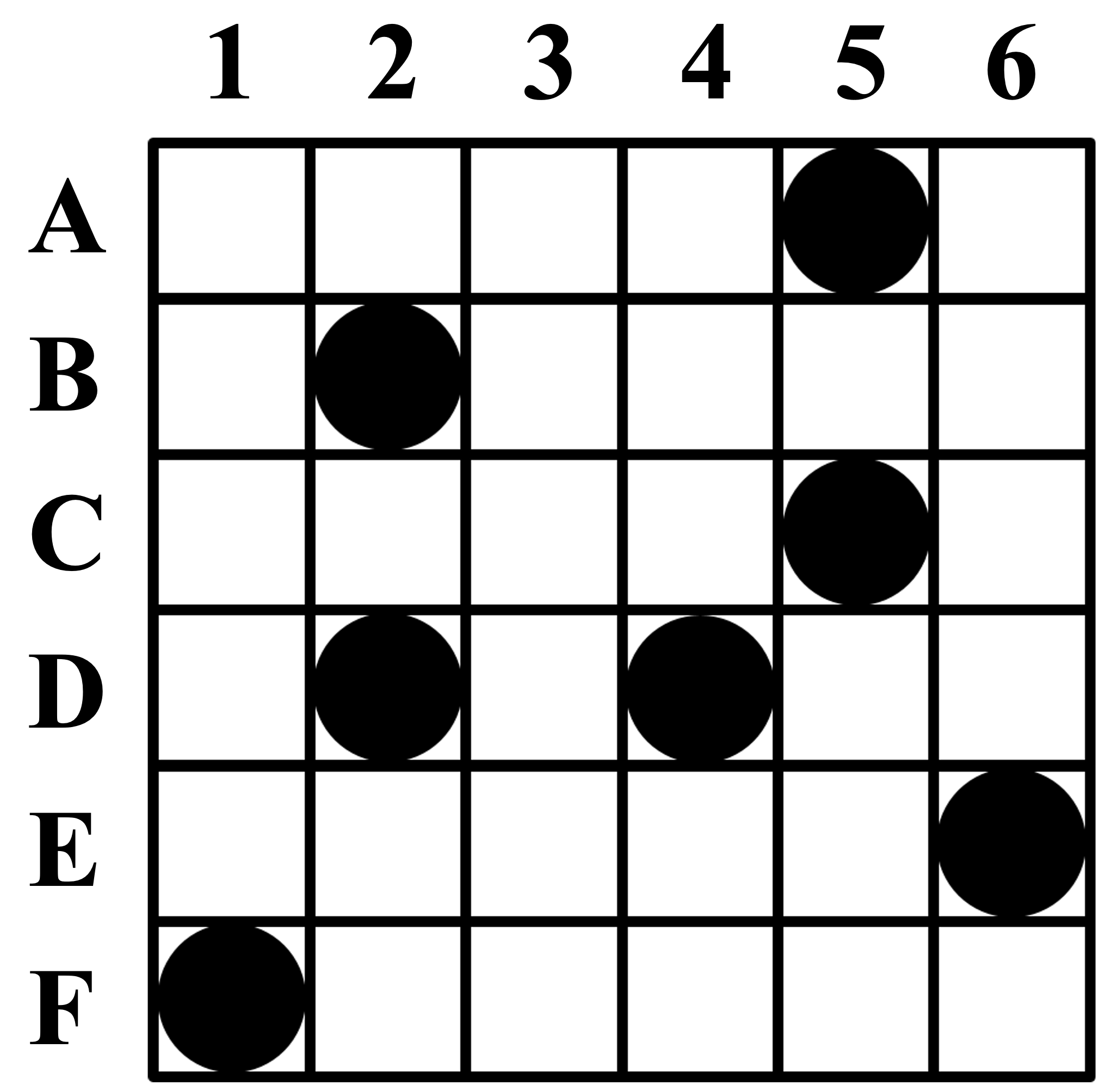}
	\end{center}
	\caption{\textit{The Genius Square} example blocker placement.}
	\label{GSEx}
\end{figure}

{\it The Genius Square} is an award-winning game published by The Happy Puzzle Company \cite{hpc}.  In the game, players create and then solve puzzles, either by playing solo or by racing other players to a solution.  This work will focus on the single-player version of this game. The puzzles are created using a $6 \times 6$ grid, a collection of seven special dice, and seven round markers that we call \emph{blockers}. Each face of each of the dice indicates a certain location on the grid (A4 or C2, for example).  When the seven dice are rolled, the player places the blockers on the seven corresponding locations, creating a puzzle that is ready to be solved. Figure \ref{GS} shows the game elements, and Figure \ref{GSEx} shows an example of blocker placement.  

\begin{figure}[h!]
	\begin{center}
		\includegraphics[width=0.3 \textwidth]{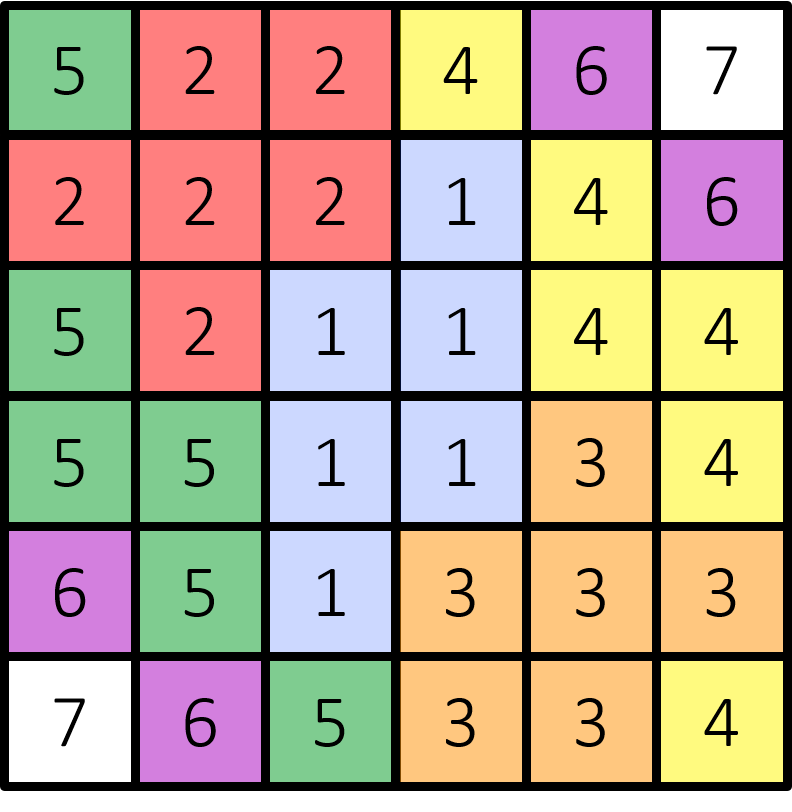}
	\end{center}
	\caption{Blocker placement options on grid from each of the seven game dice.}
	\label{Dice}
\end{figure}

We note here that each of the locations on the board is represented on exactly one of the seven dice, and no two dice have any locations in common. Figure \ref{Dice} shows the relationship between grid locations and coverage of each of the seven dice (labeled 1 through 7). Notice that die 7 only has two possible blocker locations, but most other dice have six options each. 

\begin{figure}[h!]
	\begin{center}
		\includegraphics[width=0.35 \textwidth]{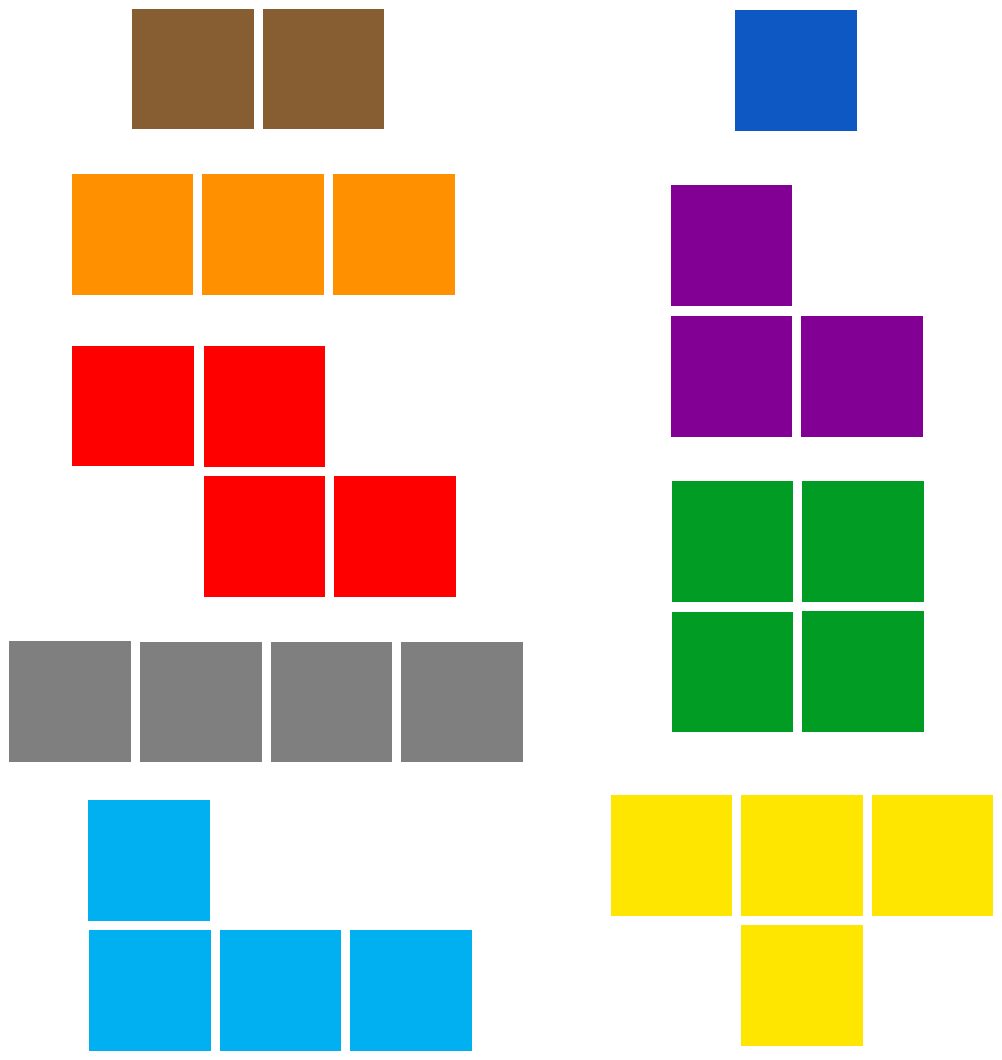}
	\end{center}
	\caption{\textit{The Genius Square} polyomino game pieces.}
	\label{Pieces}
\end{figure}

To solve the puzzle, the player tries to place nine polyomino pieces that are shown in Figure \ref{Pieces} in such a way that they cover all unblocked squares of the grid.  There are a total of 29 individual squares among all of the polyomino pieces, and, fittingly, there are 29 unblocked squares on the grid.  The game manufacturers indicate that every possible roll of the dice leads to a solvable puzzle.

%After blockers are placed on the board, the puzzler must configure the nine polyominoes to cover the remaining, non-blocked grid spaces.  Not every configuration of blockers is possible through rolling the dice, but the makers of the game claim that every configuration based on dice rolls results in a puzzle that can be solved.  An example of a blocker configuration that is not possible and not solvable is to put blockers on all of the purple grid boxes, labeled with a '6' in Figure \ref{Dice} since there is only one $1 \times 1$ polyomino piece that can be placed on the board.  
%
%The nine polyominoes are shown in Figure \ref{Pieces} and range from a $1 \times 1$ piece that occupies one grid square to pieces that take up four squares.  The total number of grid spaces contained in the polyominoes is 29.  Combined with the seven blockers, this gives a potential to cover all 36 squares of the $6 \times 6$ grid. 

%%%%%%%%%%%%%%%%%%%%%%%%%%%%%%%%%%%%%%%%%%%%%%%%%%%%%%%%%%%%%%%%%%%%%%%%%%%%%%%%%%%%%%%

\subsection{Decision Variables and Objective Function}
\label{subsec:1}
%We use an integer programming model formulation to try to assign the nine polyominoes to the non-blocked grid squares.  This is an interesting undergraduate research problem because the decision of where to place a polyomino is complicated because each piece could have several orientations.  For example, the $3 \times 1$ piece could be placed vertically or horizontally.  Further, to cover the $A2$ location on the grid, this piece could be placed so that its middle square is on this location or to the right or one space below this location.  These complications must be confronted with constraints, but the orientation problem cannot.  So, each orientation of a single polyomino must correspond to a different variable.  We further identify a location (referred to as an ``anchor dot") on each polyomino that identifies where the model will station a piece on the board.  Hence for each orientation and for each grid location there will be a binary random variable encoding whether that piece and orientation is anchored in location $(i,j)$ on the grid.  Each variable and anchor dot is seen in Figure \ref{DV}.  As an example, we define $a_1(i,j)$ and $a_2(i,j)$ to be $0-1$ variables indicating whether or not the two orientations, vertical or horizontal, of the $2 \times 1$ piece is anchored in location $(i,j)$, $1 \leq i \leq 6$, $1 \leq j \leq 6$.  Note that we have replaced the row labels $A-F$ with $1-6$ for ease of modeling.

We use an integer programming model to help assign the nine polyomino pieces to locations on the board.  To make things a little easier in terms of notation, we choose to label the board locations with (row, column) integer valued ordered pairs instead of the letter/number combinations.  The space C2, for instance, is now referred to as the space (3,2) in the third row, second column.

Another adjustment for us involves the placement of the blockers.  In the actual game, the locations of the seven blockers are completely determined by the labels on the dice (whose possible configurations are mapped out in Figure \ref{Dice}).  In our work, we will consider ALL possible placements of seven blockers, even those that would not be produced by dice rolls, which extends this research work beyond the scope of \textit{The Genius Square} game.  This means that we may encounter placements that do not allow for solutions to exist, which we discuss further in Section \ref{subsec:3}. 

There are a number of challenges for undergraduates to tackle when setting up this model.  First, the polyominoes can be placed in a number of different orientations.  For instance, the $3 \times 1$ piece (shown in orange in Figure \ref{Pieces}) could be placed either horizontally or vertically on the board.  Second, depending on a piece’s orientation, there are a certain set of locations where that piece can not be placed.

\begin{figure}[h!]
	\begin{center}
		\includegraphics[width=0.99 \textwidth]{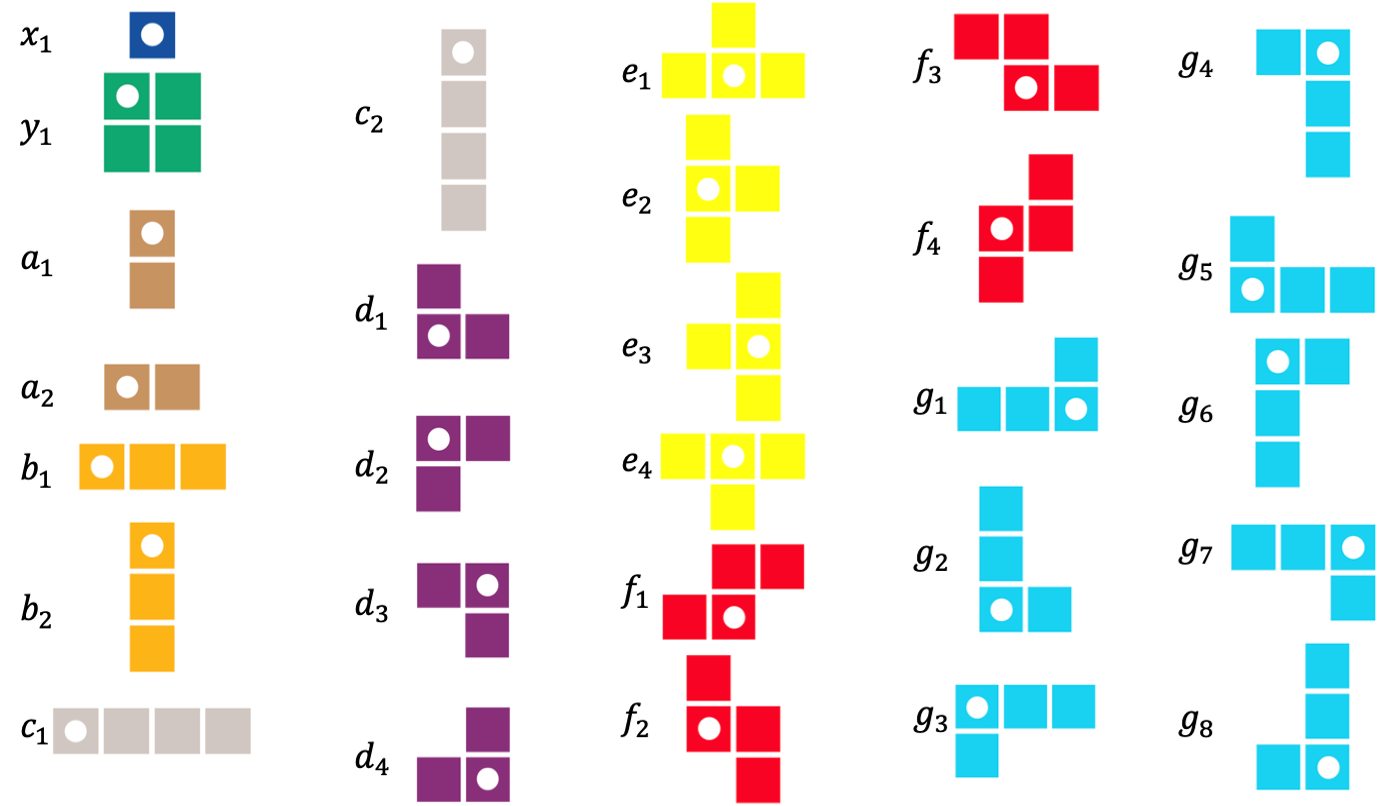}
	\end{center}
	\caption{Decision variables for each polyomino piece with anchor squares indicated by white circles.}
	\label{DV}
\end{figure}

To manage this situation, we identified a particular square on each polyomino to use as an anchor square, and this allowed us to be specific about where a piece is placed (by saying which location on the board the anchor square is located for a given piece configuration).  We created one variable for each piece/orientation/location combination.  Figure \ref{DV} shows all of the piece/orientation items as well as which of the squares is defined as the anchor square in each case (denoted by a white dot) and the decision variable name we use for each piece/orientation.  There are 28 piece/orientation pairs.  Each
such pair has an anchor dot that (in theory) could be placed in any of 36 different squares.  This gives us a total of $28 \cdot 36 = 1008$
binary variables.   These variables will be of the form $P(i,j)$ (where $P$ is one of the 28 piece/orientation pairs) and will have the value 1 if the anchor square of the piece/orientation is actually placed in position $(i,j)$, and 0 otherwise. A solution will have exactly nine of these variables with the value 1.

Before proceeding further, we  note that the nature of the piece/orientation pairs makes some dot placements impossible.  For 
instance, consider the piece/orientation labeled $a_1$ in Figure \ref{DV} (tan domino in vertical orientation). 
Since the anchor square is the top square here, it is not possible for the piece to be placed in any of the squares in row 6.   This means that six of the 36 binary variables associated with this piece/orientation will
automatically need to be 0 (and this will be handled in the constraints).

The objective function here will be set up to maximize the number of non-blocked squares that are covered.  Each polyomino in \textit{The Genius Square} game covers anywhere from one to four spaces.  This leads to an objective function of the form

%The objective function is set up to maximize coverage of the non-blocked grid.  Although, each configuration of blockers from the dice game are all solvable, we consider other configurations of blockers that are not achievable by the dice.  Hence we want to maximize our objective function, which is given by

%$$\sum_{i=1}^6 \sum_{j=1}^6 \Large(x_1(i,j) + y_1(i,j) + \sum_{k=1}^2 a_k(i,j) + \sum_{k=1}^2 b_k(i,j) + \\
%\sum_{k=1}^2 c_k(i,j) +\sum_{k=1}^4 d_k(i,j) +\sum_{k=1}^4 e_k(i,j) + \sum_{k=1}^4 f_k(i,j) + \sum_{k=1}^8 g_k(i,j) \Large).$$
%
%\christian{Below is my attempt to write the same objective function but make it more compact.}

%$$\sum_{i=1}^6 \sum_{j=1}^6 \left(x_1(i,j) + y_1(i,j) + \sum_{k=1}^2 \left(a_k(i,j) +  b_k(i,j) + c_k(i,j)\right) +\sum_{k=1}^4 (d_k(i,j) + e_k(i,j) + f_k(i,j)) + \sum_{k=1}^8 g_k(i,j) \right).$$

\[
\begin{array}{ll}
\displaystyle \sum_{i=1}^6 \sum_{j=1}^6 & \displaystyle \bigg(x_1(i,j) + 4y_1(i,j) + \sum_{k=1}^2 \left[2a_k(i,j) +  3b_k(i,j) + 4c_k(i,j)\right]  \\ & \displaystyle  + \sum_{k=1}^4 \left[ 3d_k(i,j) + 4e_k(i,j) + 4f_k(i,j) \right] + \sum_{k=1}^84 g_k(i,j) \bigg),
\end{array}
\]
which we want to maximize. Note that this objective function will allow for solutions that successfully cover the board as well as ones that are unsolvable, which will be discussed further in Section \ref{subsec:3}.

%%%%%%%%%%%%%%%%%%%%%%%%%%%%%%%%%%%%%%%%%%%%%%%%%%%%%%%%%%%%%%%%%%%%%%%%%%%%%%%%%%%%%%

\subsection{Constraints}
\label{subsec:2}

Constraints for \emph{The Genius Square} problem involve ensuring no piece hangs off the board, covers another piece, or covers a blocker.  We will not cover every type of constraint here, but we offer example constraints to give the reader an idea. We leave other constraints as exercises.  

Consider the piece whose decision variables are labeled $d_k(i,j)$ for $k=1,2,3,4$ (purple in Figure \ref{DV}).  Since the anchor square for this piece is the center piece, it is possible to anchor this piece in a certain orientation where one square of this piece would hang off the board.  For example, if $d_3(2,1)=1$ then one square of this polyomino hangs off the left side of the grid. %does not cover a grid square when position in location $(2,1)$.  
For each piece and each orientation where we could have this situation occur, we establish a constraint that avoids the situation.

\begin{exercise}
Write a constraint for the $d_3(i,j)$ variables, $1 \leq i,j \leq 6$ pieces indicating that this piece/orientation can be used at most once on the board. 
\end{exercise}

For each non-blocked square, we ensure that two different orientations orientations of the same piece are not both played.  For example, if location $(1,1)$ is not occupied by a blocker, we set up the following constraint:

%$$x_1(1,1) +y_1(1,1) + \sum_{k=1}^2 a_k(1,1) +\sum_{k=1}^2 b_k(1,1) + \sum_{k=1}^2 c_k(1,1) + \\
%d_1(2,1) + d_2(1,1) + d_3(1,2) +e_2(2,1) +e_4(1,2) + f_2(2,1) + f_3(2,2) +g_2(3,1) + g_3(1,1) + g_4(1,2) +g_5(2,1) +g_6(1,1) +g_7(1,3) \leq 1.$$

\[
\begin{array}{ll}
\displaystyle \bigg( x_1(1,1) +y_1(1,1) + \sum_{k=1}^2 \left[ a_k(1,1) + b_k(1,1) + c_k(1,1)  \right]&\\
\displaystyle \mbox{\ \ \ \ }+ d_1(2,1) + d_2(1,1) + d_3(1,2) + e_2(2,1) +e_4(1,2) + f_2(2,1)  + f_3(2,2) &\\ 
\displaystyle \mbox{\ \ \ \ } +g_2(3,1) + g_3(1,1) + g_4(1,2) +g_5(2,1) +g_6(1,1) +g_7(1,3) \bigg)&\leq 1.
\end{array}
\]
This type of constraint is repeated for each non-blocked space. By using an inequality instead of an equation equal to 1, we are allowing for the possibility that a board is not solvable, which would mean that there are unfilled spaces and unused pieces. For blocked spaces, we produce a similar constraint except we restrict the sum to equal zero.  Finally, we ensure that only one configuration of each piece is played at most once since there are no duplicate pieces provided with \textit{The Genius Square} game.  For example, for the blue piece labeled with decision variable $g_k(i,j)$, this constraint looks like

$$\sum_{i=1}^6 \sum_{j=1}^6 \sum_{k=1}^8 g_k(i,j) \leq 1.$$
 Note that additional constraints keep this from causing an infeasibility, such as setting $g_3(6,6)=1$.

\begin{exercise}
Write a constraint that ensures the piece $x_1$ (single square) must be used exactly once.  The decision variables for this piece are $x_1(i,j)$ for $1 \leq i,j \leq 6$.
\end{exercise}

\begin{exercise}
Write a constraint that ensures the domino piece $a$ (two squares) must be used exactly once.  The decision variables for this piece are $a_1(i,j)$ and $a_2(i,j)$.  Be careful not to place this piece in a location where one square hangs off the board and make sure to specify allowable values of $i$ and $j$.  

\end{exercise}

%%%%%%%%%%%%%%%%%%%%%%%%%%%%%%%%%%%%%%%%%%%%%%%%%%%%%%%%%%%%%%%%%%%%%%%%%%%%%%%%%%%%%%%%%%%%

\subsection{Results}
\label{subsec:3}
\begin{figure}
     \centering
     \begin{subfigure}[b]{0.3\textwidth}
         \centering
         \includegraphics[width=\textwidth]{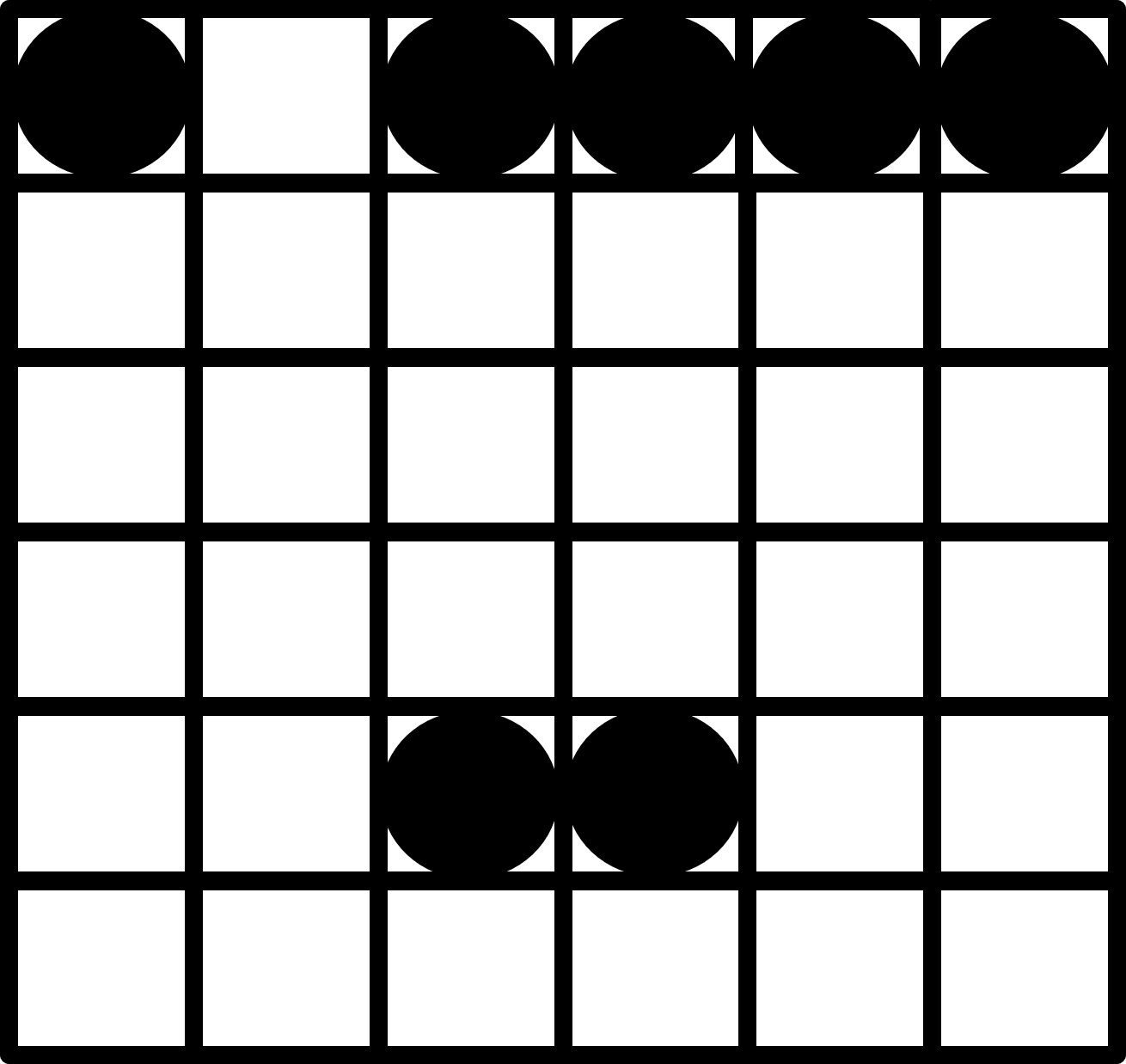}
	\caption{22,317 solutions.}
	\label{MS}
     \end{subfigure}
   \mbox{\ \ \ \ \ \ \ \ }
     \begin{subfigure}[b]{0.3\textwidth}
         \centering
         \includegraphics[width=\textwidth]{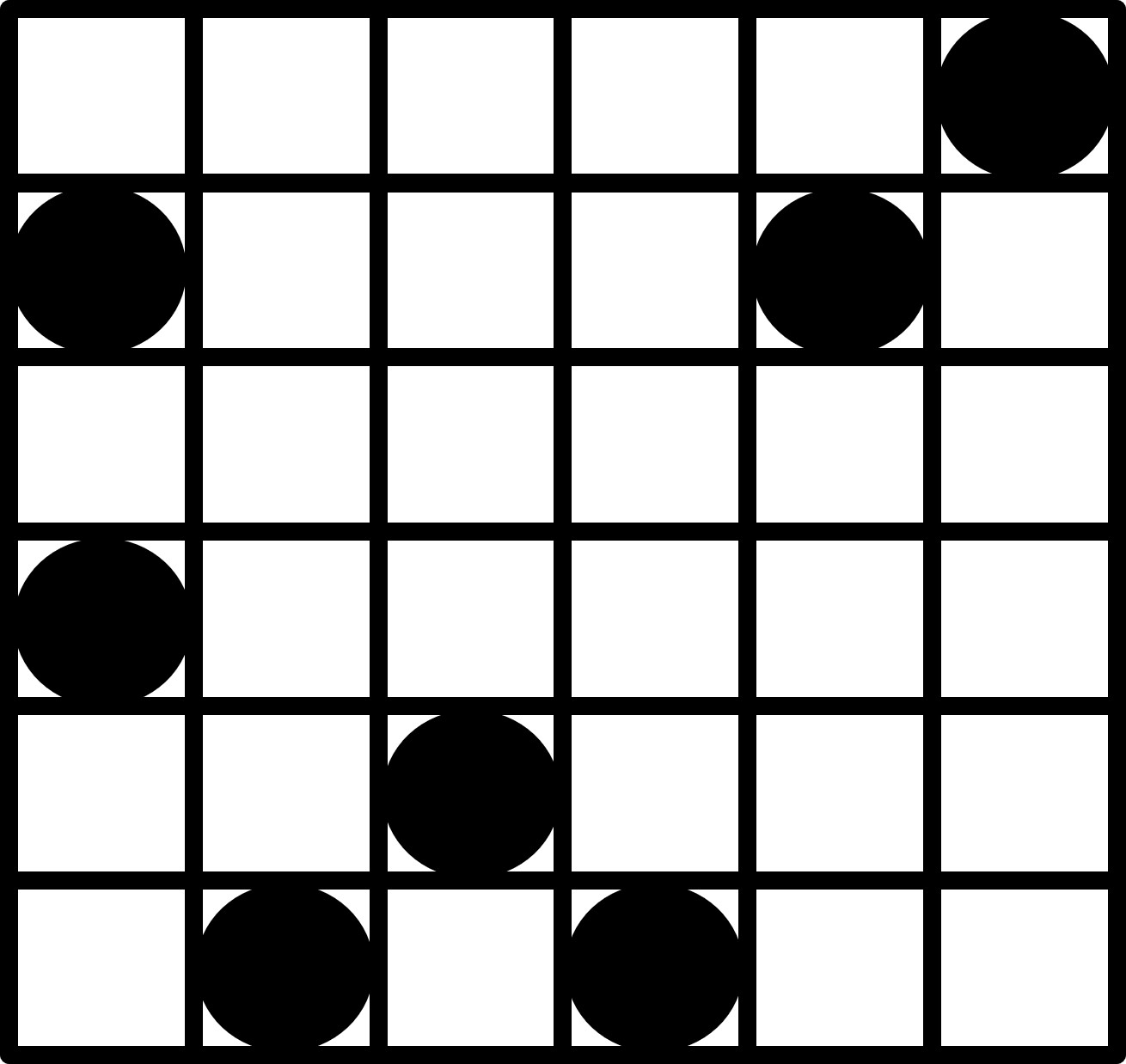}
	\caption{11 solutions.}
	\label{FS}
     \end{subfigure}
        \caption{Game configurations with (a) many solutions and (b) the minimum number of solutions (for game boards generated by the game's dice).}
        \label{fig:Solutions}
\end{figure}

Using this model and the Gurobi Optimization software package, our students confirmed that the manufacturer's claim that every blocker configuration generated by the dice leads to a solution.  In fact most games have many solutions, sometimes tens of thousands.  For example, the blocker configuration shown in Figure \ref{MS} has 22,317 solutions, while some solvable puzzles have very few solutions. Figure \ref{FS} shows an example of a blocker configuration with only 11 solutions.

%\begin{figure}[h!]
%	\begin{center}
%		\includegraphics[width=0.25 \textwidth]{ManySolutions.jpg}
%	\end{center}
%	\caption{A game configuration with over 22,000 solutions}
%	\label{MS}
%\end{figure}
%
%\begin{figure}[h!]
%	\begin{center}
%		\includegraphics[width=0.25 \textwidth]{FewSolutions.jpg}
%	\end{center}
%	\caption{A game configuration with 11 solutions}
%	\label{FS}
%\end{figure}

%We'll do an example with three images along side each other with separate captions and labels. Here's some example code:

\begin{figure}
     \centering
     \begin{subfigure}[b]{0.3\textwidth}
         \centering
         \includegraphics[width=\textwidth]{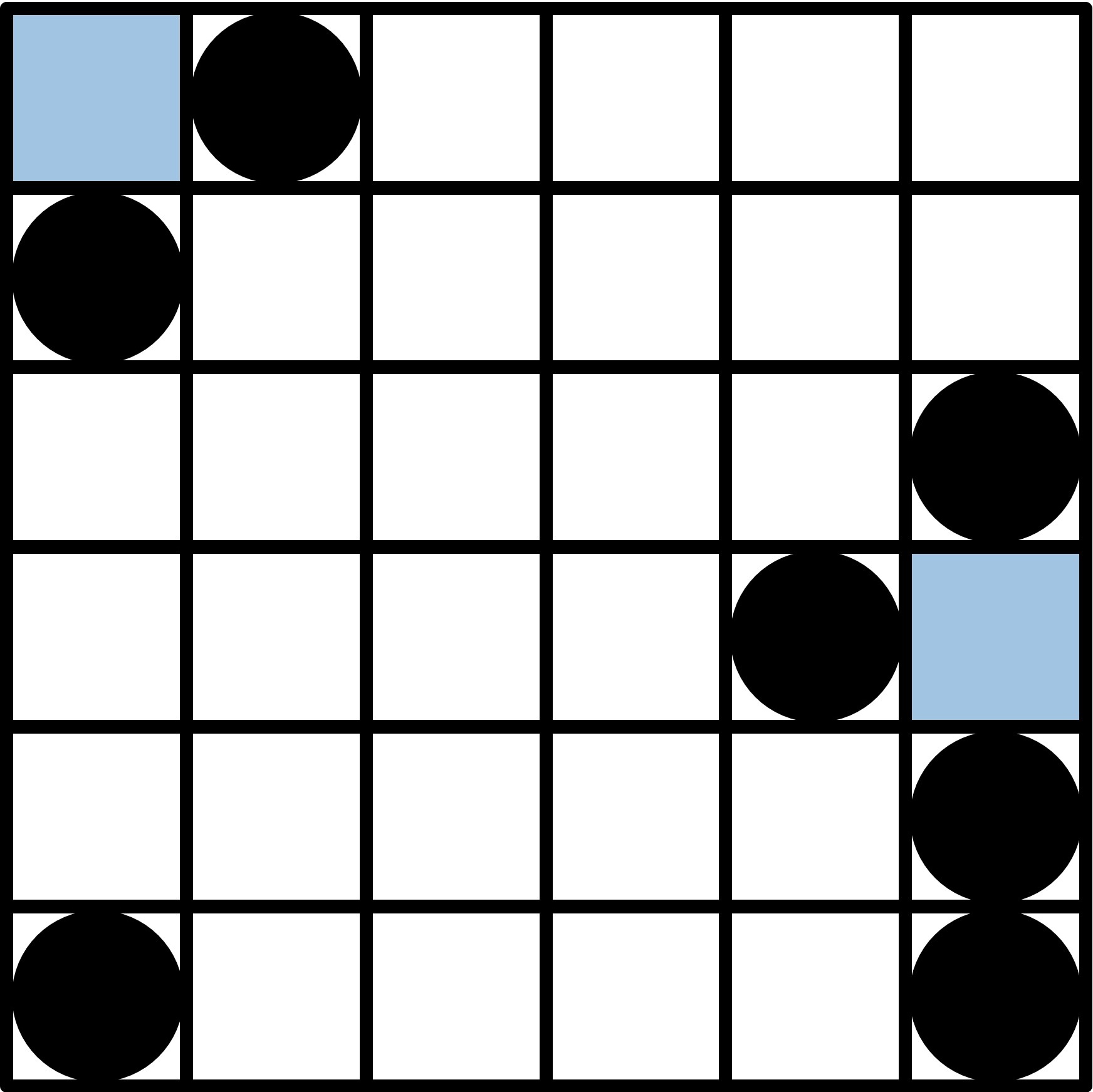}
	\caption{Constrained one-piece scenario}
	\label{FI1}
     \end{subfigure}
   \mbox{\ \ \ \ \ \ \ \ }
     \begin{subfigure}[b]{0.3\textwidth}
         \centering
         \includegraphics[width=\textwidth]{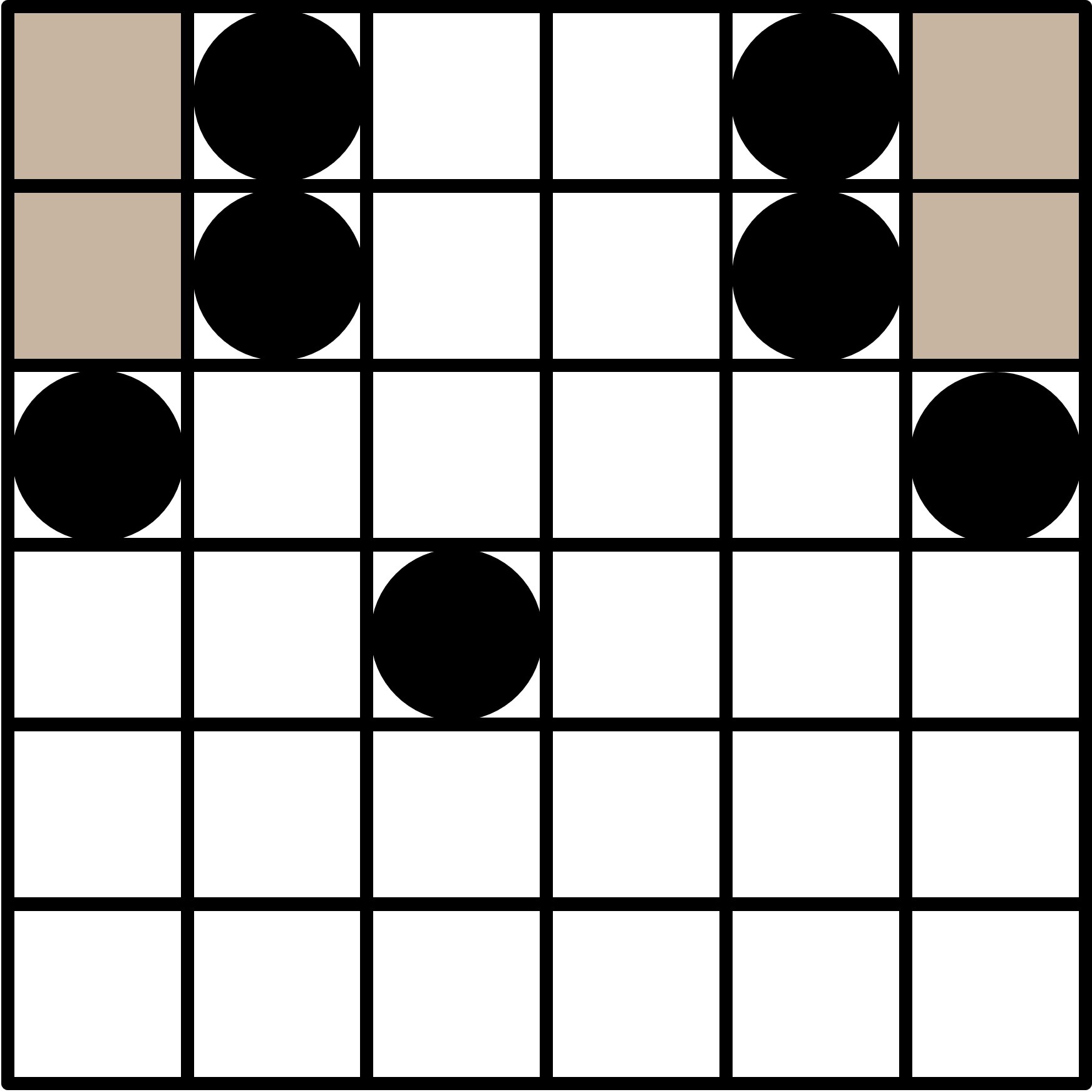}
	\caption{Constrained two-piece scenario}
	\label{FI2}
     \end{subfigure}
        \caption{Examples of non-solvable games because player is forced to use the same piece twice.}
        \label{fig:FI}
\end{figure}

As mentioned earlier, we tested blocker configurations that are not achievable with the game dice as part of our research work. We found that most of these configurations have solvable game puzzles as well.  Of the over 1.04 million puzzles, approximately $2\%$ of those are unsolvable.  Most of the unsolvable boards can be explained by one of three scenarios.  The first scenario, which we will call the \emph{constrained one-piece scenario}, creates a blocker configuration that puts blockers in locations that force the player to use the $1 \times 1$ piece twice.  For example, consider the blocker placement shown in Figure \ref{FI1}.
%placing blockers in locations $(1,5), (2,6), (5, 1),$ and $(6,2)$.  
This forces two $1 \times 1$ game pieces to be played in the squares shaded blue, but the game only includes one of these pieces.  In the second scenario shown in Figure \ref{FI2}, which we will call the \emph{constrained two-piece scenario}, the puzzler is forced to play two $2 \times 1$ game pieces, which is impossible because only one is included with \emph{The Genius Square} game.  %Figure \ref{FI} shows an example of the Constrained Two-Piece Scenario and a different example of a Constrained One-Piece Scenario.

%\begin{figure}[h!]
%	\begin{center}
%		\includegraphics[width=0.45 \textwidth]{ForcedImpossibilities.jpg}
%	\end{center}
%	\caption{Example of non-solvable game because player is forced to use the same piece twice}
%	\label{FI}
%\end{figure}

\begin{figure}
     \centering
     \begin{subfigure}[b]{0.3\textwidth}
         \centering
         \includegraphics[width=\textwidth]{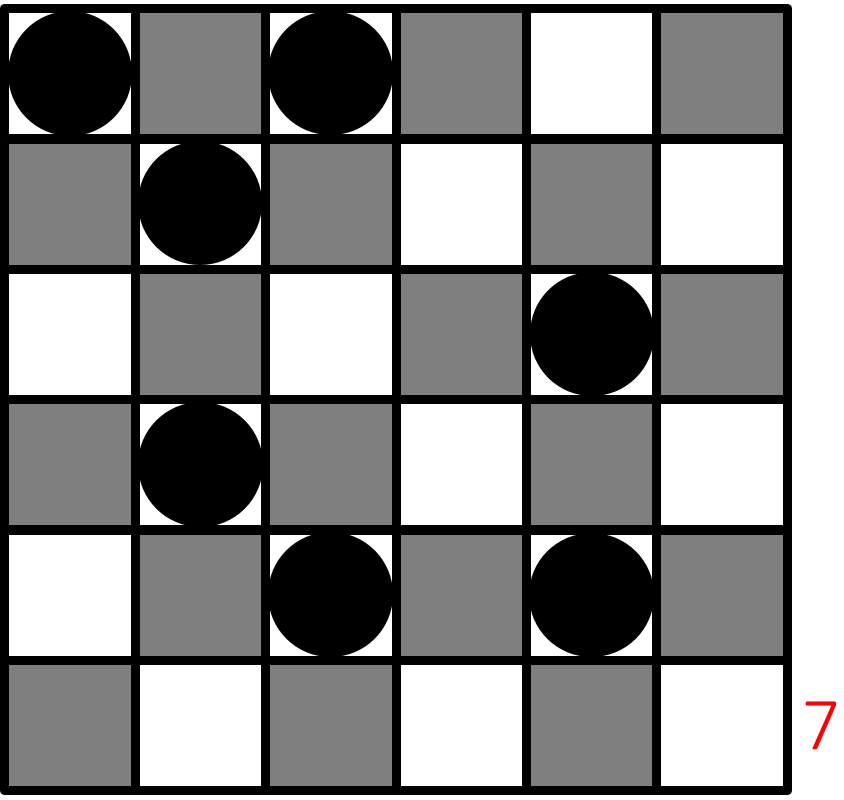}
	\caption{Placing all blockers on white squares leaves seven more gray squares than white squares.}
	\label{PIBoard}
     \end{subfigure}
   \mbox{\ \ \ \ \ \ \ \ }
     \begin{subfigure}[b]{0.3\textwidth}
         \centering
         \includegraphics[width=\textwidth]{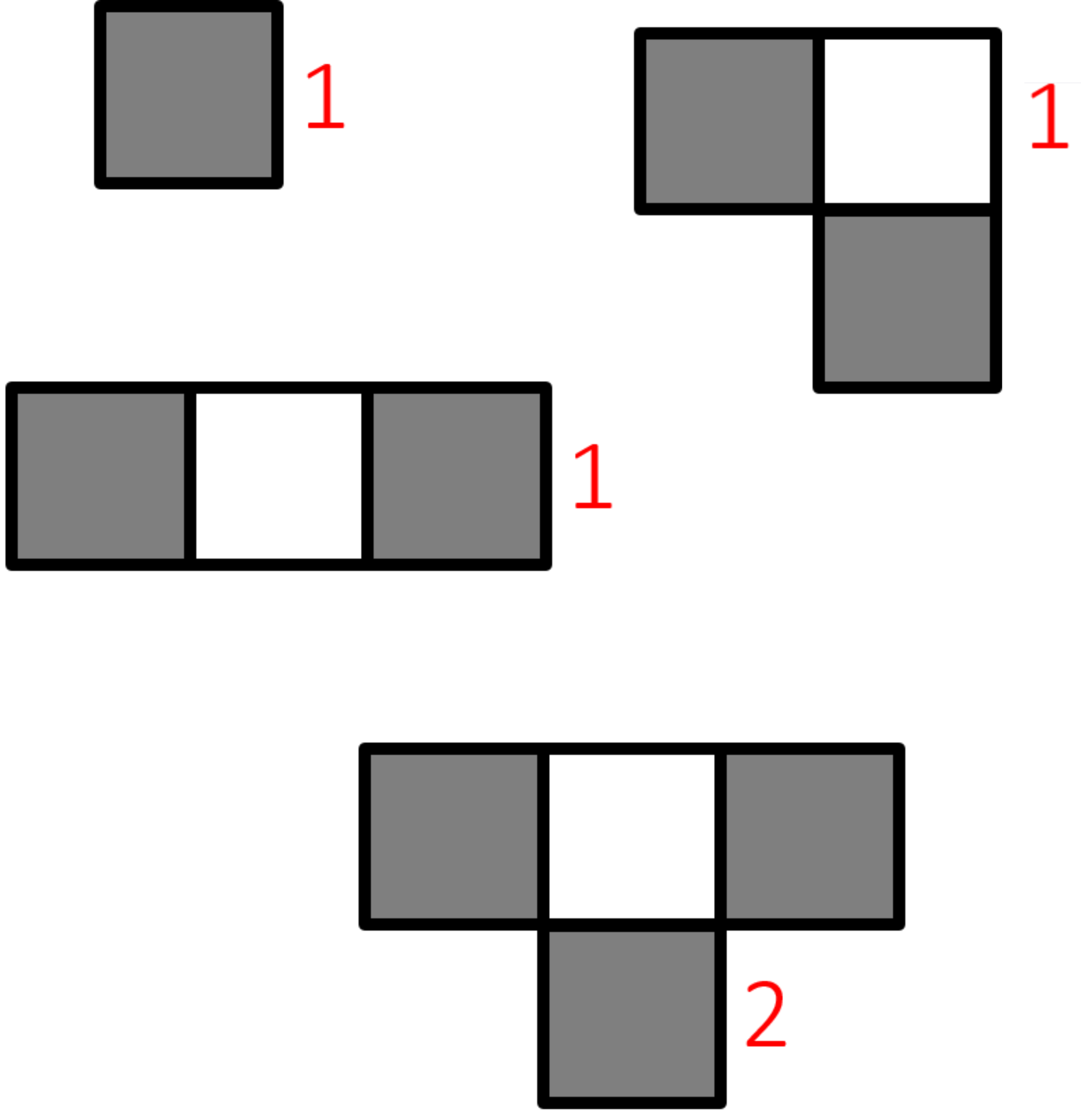}
	\caption{Of the polyominoes with an imbalance of alternating colored squares, the imbalance can be at most five. All pieces not pictured have the same number of white and gray squares. }
	\label{PIPieces}
     \end{subfigure}
        \caption{Examples of parity imbalance scenario where the game board has an imbalance of seven, but the pieces have an imbalance of five.}
        \label{fig:PI}
\end{figure}

\begin{figure}
     \centering
     \begin{subfigure}[b]{0.3\textwidth}
         \centering
         \includegraphics[width=\textwidth]{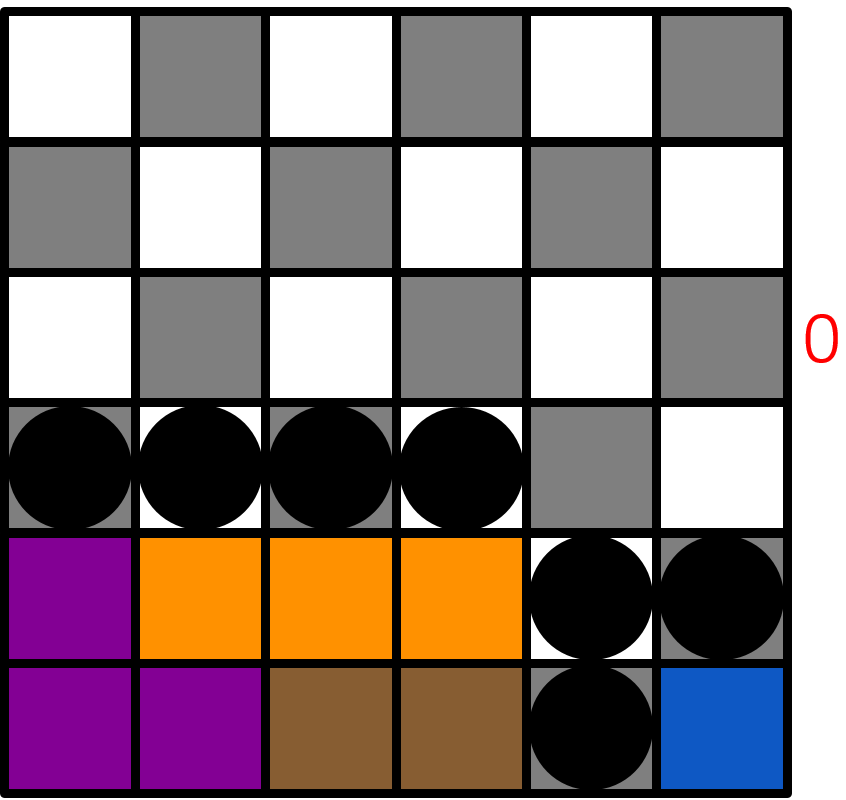}
	\caption{Game board with balance among remaining gray and white squares}
	\label{ComplexPIBoard}
     \end{subfigure}
   \mbox{\ \ \ \ \ \ \ \ }
     \begin{subfigure}[b]{0.3\textwidth}
         \centering
         \includegraphics[width=\textwidth]{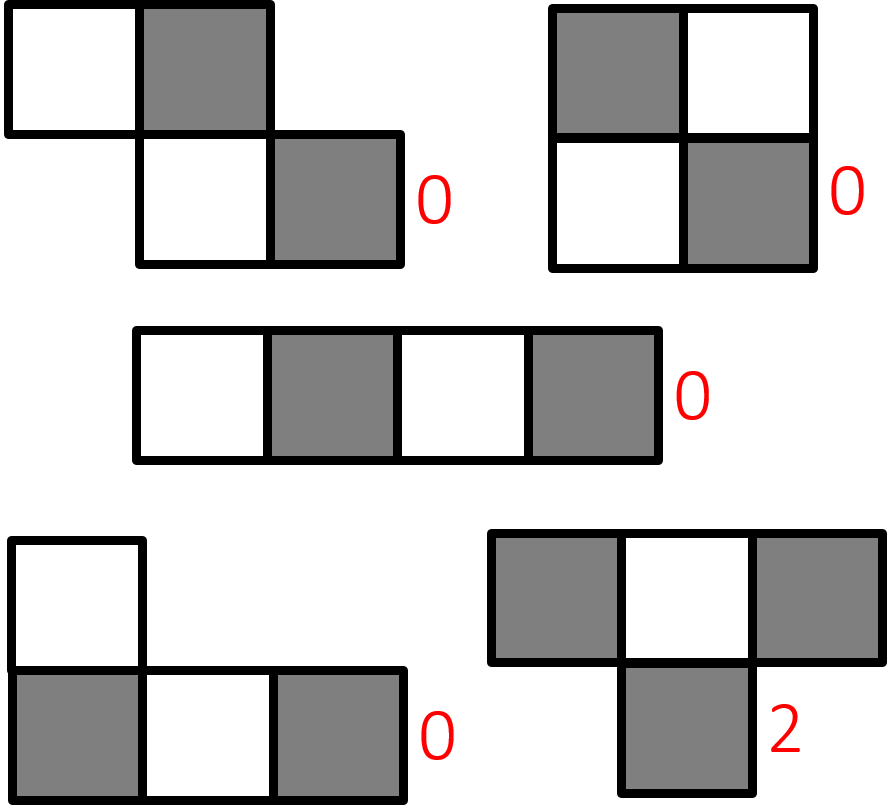}
	\caption{Remaining pieces have a net parity of two}
	\label{ComplexPIPieces}
     \end{subfigure}
        \caption{A parity argument where the difference among remaining pieces is two, but the board has parity balance.}
        \label{fig:ComplexPI}
\end{figure}

Most of the other impossible-to-solve game boards are from a blocker configuration that creates a parity imbalance.  We will call this the \emph{parity imbalance scenario}.  To see this, envision a checkerboard configuration of the $6 \times 6$ game board as shown by the gray and white shading in Figure \ref{PIBoard}.  Some of the polyominoes, like the $2 \times 2$ square, cover the same number of black and white checkerboard squares.  However, the four polyominoes shown in Figure \ref{PIPieces} cover an unequal number of gray and white squares.  So, if a blocker configuration, such as the one shown in Figure \ref{PIBoard}, is placed so that it covers squares all of the same color, that leaves a deficit of the other color's squares that have to be recouped using the game pieces.  However, the four game pieces where there is an imbalance in the colors can only recoup five of that color, creating an infeasible configuration since the game board has an imbalance of seven.  
%
%\begin{figure}[h!]
%	\begin{center}
%		\includegraphics[width=0.45 \textwidth]{ParityImbalanceExample.png}
%	\end{center}
%	\caption{Example of Parity Imbalance}
%	\label{PI}
%\end{figure}
%
Figure \ref{fig:ComplexPI} shows a more complex scenario where there are not enough imbalanced game pieces left to recoup the color imbalance.

%\begin{figure}[h!]
%	\begin{center}
%		\includegraphics[width=0.6 \textwidth]{ComplexParityImbalance.png}
%	\end{center}
%	\caption{Complex scenario of parity imbalance}
%	\label{CPI}
%\end{figure}

\begin{figure}[h!]
	\begin{center}
		\includegraphics[width=0.3 \textwidth]{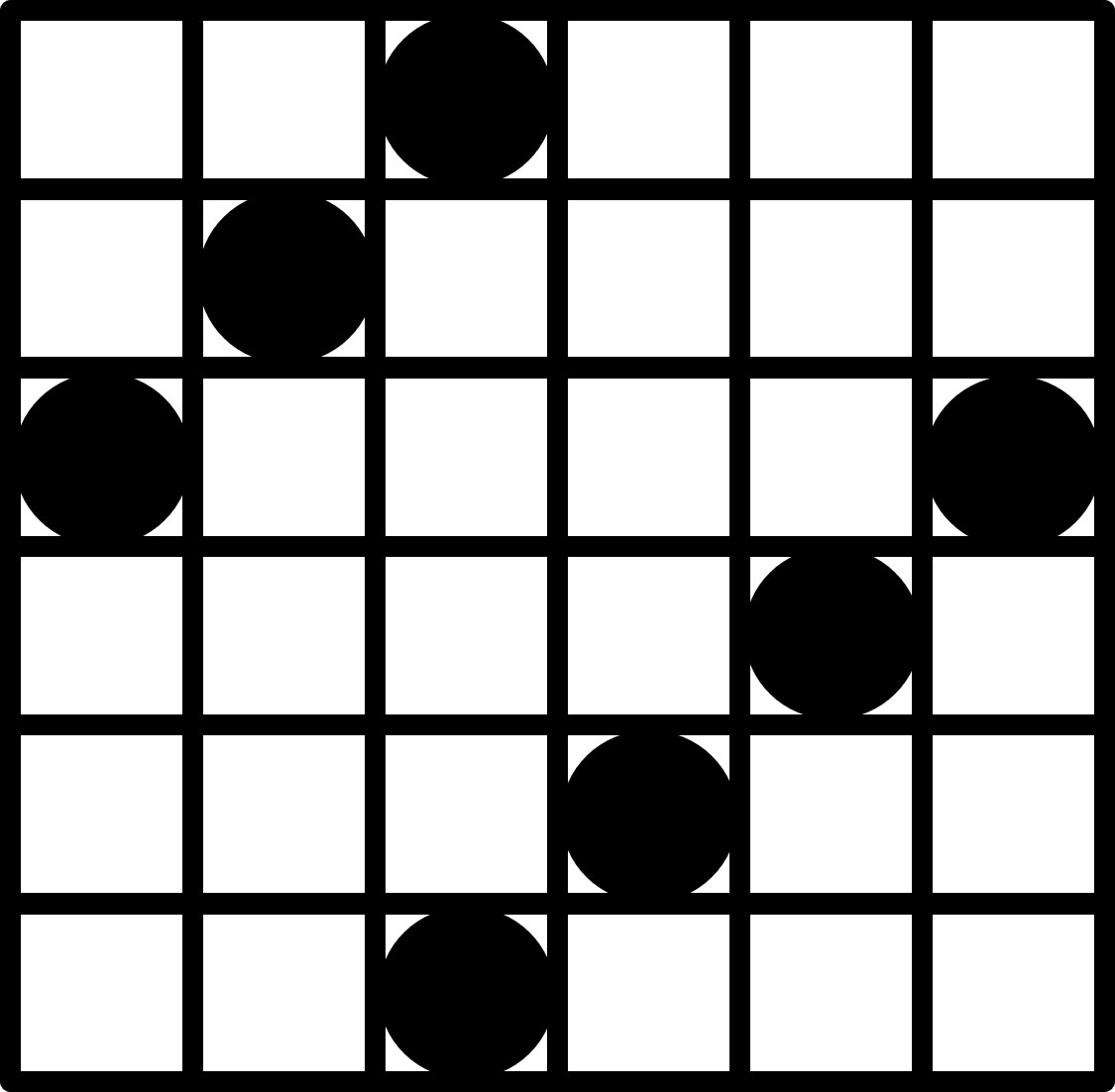}
	\end{center}
	\caption{An unsolvable board that doesn't fall into a constrained piece scenario or a parity imbalance scenario}
	\label{CFP}
\end{figure}

Even with that, there are still 148 game configurations that create an unsolvable game that can't be explained by one of the three aforementioned scenarios.  This leaves an open area for more exploration.  Figure \ref{CFP} shows an example of a game configuration that does not fall into the scenarios listed above.  Explaining why this type of configuration does not lend itself to a solution is an open problem.  

\resproject{
Come up with proof of the unsolvability of the complex scenarios.}

\resproject{\emph{The Genius Star} is a variant of \emph{The Genius Square} that has a star configuration for the game board.  The game play is similar in that the player rolls dice, places blockers, and then tries to place all of the game pieces onto the unblocked locations  This game can be analyzed similarly to find (i) how many solutions are possible for each blocker configuration and (ii) whether there exist blocker configurations where solutions are not possible.}

%%%%%%%%%%%%%%%%%%%%%%%%%%%%%%%%%%%%%%%%%%%%%%%%%%%%%%%%%%%%%%%%%%%%%%%%%%%%%%%%%%%%%%%%%%%%%%

\section{Ticket to Ride}
\label{sec:4}
\begin{figure}[h!]
\begin{center}
\includegraphics[width = 1.0\textwidth]{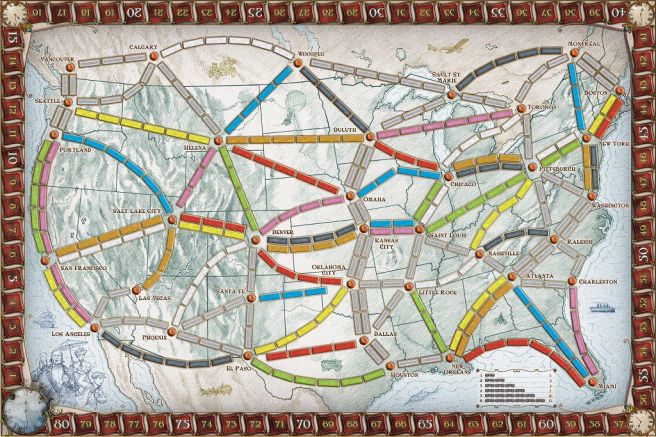}
\caption{\emph{Ticket to Ride} game board.  Image courtesy of  \cite{bgg}.}
\label{TTR}
\end{center}
\end{figure}

% Always give a unique label
% and use \ref{<label>} for cross-references
% and \cite{<label>} for bibliographic references
% use \sectionmark{}
% to alter or adjust the section heading in the running head
\emph{Ticket to Ride} is a popular board game published by Days of Wonder \cite{dow}. In the game, players compete against each other in allocating train resources to routes in order to maximize the number of points they receive in the game.  Each player is given 45 trains that they can use to claim routes on the game board spaces seen in Figure \ref{TTR}.  To allocate trains on the board, a player must accumulate cards of the same number and color, or some combination of that color and wild cards, as the desired route of the game board. For example, to claim the route from Miami to Charleston requires the player to have four pink cards (or a combination of pink and wild cards) and four trains to claim the route.  For each route claimed, a player gets a certain number of points in the game based on the length of the route as indicated in Table \ref{C25}.

\begin{table}[htbp!] 
\begin{center}
\begin{tabular}{|c|c|} \hline
 Route Length & Points    \\ \hline
 1 & 1  \\ \hline
 2 & 2 \\ \hline
 3 & 4  \\ \hline
 4 & 7 \\ \hline
 5 & 10 \\ \hline
 6 & 15 \\ \hline
\end{tabular}
\caption{Points per route length.}
\label{C25}
\end{center}
\end{table}

A player also has at least two destination cards in the game that have point values assigned to them.  If the owner of the destination card creates a path in the game with their trains between the destination cities on the card by the end of the game, the player gains the amount of points indicated on the card.  Otherwise the player loses those same number of points.  For example, one destination card has the cities Seattle to New York and scores 22 points if a path between the cities is completed by the owner of this card by the end of the game and loses 22 points if not.

This game has been studied by several authors from the graph theoretic point of view, some employing Dijkstra's Algorithm to find shortest paths through the game board network and breaking down points per route (e.g., \cite{Bettles, Reiber, Slade, WitterLyford}).  Our research goal was different and two-fold.  We wanted to set up an integer program to determine the optimal allocation of trains to routes to maximize the number of points in the game connecting all of the destinations on a given set of three destination cards.   We make the simplifying assumption that all trains played go towards the purpose of connecting the destination cities, although in actual game play, players might choose to forego connecting their destination cities in order to claim edges that award higher points.  We then wanted to iterate our analysis over any set of three or four destination cards.  In each case, we sought to find the routes that were used most often in forming paths in those optimal solutions.  It's worth noting that we were not replicating game play scenarios because we were working under the assumption that no strategic game play from one's opponents can impact our train placement optimization, where in reality, opponents' trains may make certain routes unavailable. 

This problem is well-suited for an undergraduate research project.  In addition to having students investigate the concepts around integer programming, the game board forms a network.  Modeling the formation of paths between destination cities is complicated by the condition that there must be constraints to ensure the individual routes (edges in the network) connect to form a path between destination cities.  Often in network optimization, edges are directed from one city to the next so that in forming a path between cities we talk about ``entering" and ``leaving" a node (city) via one-way directed edges.  However, in the \emph{Ticket to Ride} game there is often no obvious direction between cities, and we consider the edges undirected.  Thus, in trying to connect two destination cities via a path, if a city is a non-destination city, it must have either zero or two edges chosen around it.  If zero edges are chosen, then the city is not on the path between two destination cities.  If two edges are chosen, then the city is on the path.  Viewed another way, we can make a choice of whether to include a city as part of a path to connect destination cities.  If the choice is to include the city, then we choose two neighboring edges to be on the path.  If the choice is to not include the city, then we choose zero edges.  For destination cities, we must choose one edge neighboring each of these locations to include in the path.

\begin{figure}[h!]
\begin{center}
\includegraphics[width = 0.3\textwidth]{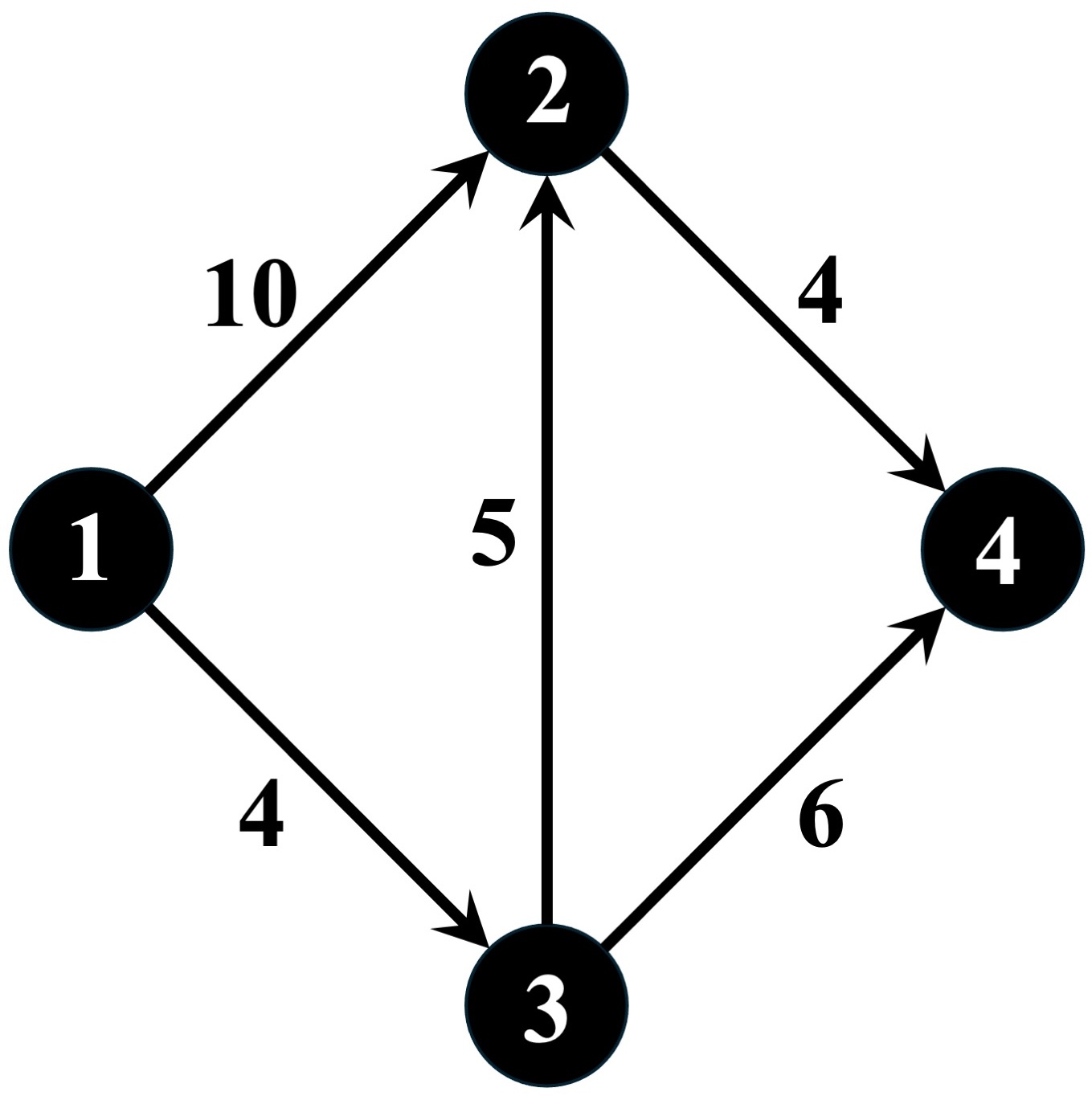}
\caption{Path example network.}
\label{sp1}
\end{center}
\end{figure}

\chproblem{Write an IP model to find a path from location one to location four in the network shown in Figure \ref{sp1} that maximizes points.  The weights (numbers) on each edge are the point values given for the edge.   Be sure to include variables for making the choice of including node two or node three in the created path.  Also include constraints to ensure that the answer is actually a path and not just a collection of the highest value edges.}

%%%%%%%%%%%%%%%%%%%%%%%%%%%%%%%%%%%%%%%%%%%%%%%%%%%%%%%%%%%%%%%%%%%%%%%%%%%%%%%%%%

\subsection{Decision Variables and Objective Function}
\label{subseq:4}

We assume that there is some arbitrary labeling of cities given.  We will describe the scenario where a player has three destination cards.  Other scenarios are similar.  With three destination cards, which we will label $A$, $B$, and $C$, a player must connect the three different pairs of cities via some set of paths of their trains.  Potentially these paths could have overlap, i.e. we could allocate trains to the network edge between New Orleans and Atlanta to be used on the paths connecting the cities on destination card $A$ and the cities on destination card $B$.  For this reason, we differentiated decision variables based on what route, either $A$, $B$, or $C$, was being covered.  However, differentiating these decision variables by route might lead to double counting the points for an edge or double counting the total number of trains used.  To get around this we set up an overall decision variable for each network edge to be used in the objective function and in the train inventory constraint.  We will set up constraints below that ``turn on" (set value to 1) this overall decision variable if at least one of the specific route variables uses the route.  

To this end, let $a_{ij}, b_{ij},$ and $c_{ij}$ be binary decision variables for each possible city $i$ and city $j$ adjacent in the train network.  A value of $a_{ij}=1$ indicates that edge $ij$ is part of the path used to connect destination cities on card $A$, and $a_{ij}=0$ otherwise (and similarly for $b_{ij}$ and $c_{ij}$).  Further, we set up variable $x_{ij}$ to be a binary decision variable for each possible adjacent cities $i$ and $j$.  We will constrain $x_{ij}=1$ if any of $a_{ij}, b_{ij}$, or $c_{ij}$ is assigned a value of 1 to identify whether a city pair shows up on some destination card in play.  Finally, we set up a variable $y^a_{i}$ (likewise $y^b_i$ and $y^c_i$) to be a binary variable representing whether city $i$ is used along path $A$ or not.  This variable in some sense will ``open" a city to include a path through it.

%\john{I think some examples would be helpful here for the reader to get a better sense of what these variables represent.}

Given a set of three destination cards, our objective is to allocate 45 trains to the board to maximize the number of points scored.  This is an ideal scenario that assumes that there are no other players claiming routes at the same time.  In other words, our solution would represent what would happen if we played optimally without interference and could achieve all the color cards that we would need.  If we let $p_{ij}$ be a parameter representing the amount of points a player receives from claiming route $ij$, then the objective function for our model is represented by

$$\sum_i \sum_j p_{ij} x_{ij},$$
which we want to maximize.  

%%%%%%%%%%%%%%%%%%%%%%%%%%%%%%%%%%%%%%%%%%%%%%%%%%%%%%%%%%%%%%%%%%%%%%%%%%%%%%%%%%%%%%%%%%%%%%%

\subsection{Constraints}
\label{subseq:5}
We will describe the constraints for enforcing a path connecting cities on destination card $A$ in detail.  The constraints for destination cards $B$ and $C$ are similar and will be omitted.  The first constraint ensures that we connect a path to our destination cities, which we will call $\alpha$ and $\beta$:
\[ 
\begin{array}{l}
\displaystyle \sum_i (a_{\alpha i} + a_{i \alpha}) = 1, \mbox{ and}\\
\displaystyle \sum_j (a_{\beta j} + a_{j \beta}) = 1.
\end{array}
\]
For all other cities, we enforce that the number of path edges used is either zero or two, depending on whether we choose to ``open" the city or not:

$$\sum_i (a_{ki} + a_{ik}) = 2y^a_i, \forall k \neq \alpha, \beta.$$

We need to ensure that we don't use more trains than we have (45).  As mentioned before, it is possible for each of the paths connecting destination cards $A$, $B$, and $C$ to use the same route.  For instance, it is possible that $a_{ij} = b_{ij}=1$ in an optimal solution.  With this mind, it is not as simple as summing over all variables $a_{ij}, b_{ij}, c_{ij}$ to total the number of trains used due to double-counting.  So, we use $x_{ij}$ for this total but need to connect $x_{ij}$'s value to be one only if at least one of $a_{ij}, b_{ij}$, or $c_{ij}$ is equal to one.  Given that $x_{ij}$ is a binary variable, we can use the following constraint for this purpose:

$$x_{ij} \geq \frac{a_{ij} + b_{ij} +c _{ij}}{3}.$$
Note that from the definitions of the variables, if $a_{ij}=b_{ij}=c_{ij}=0$, $x_{ij}$ does not necessarily have to be 0.   If a solution can be found the connects all of the destination cities, and doesn't use the route $ij$, the solver might still choose to lay trains on this edge in an optimal solution if there are trains remaining after connecting the cites.

If we let $t_{ij}$ be a parameter (not a variable) representing the number of trains needed to claim route $ij$, the constraint that we can't use more than our game allotment of trains is then

$$\sum_i \sum_j t_{ij} x_{ij} \leq 45.$$
Finally, we don't want the model to use the same route in a forward $ij$ and backward $ji$ direction because this might accumulate the same points twice which is not allowed.  
\begin{exercise}
Using the network in Figure \ref{sp1}, describe a constraint that will avoid using the same path forwards and backwards.
%The following equation holds for each route $ij$:
%
%$$x_{ij} + x_{ji} \leq 1.$$

\end{exercise}
%Once this model is formulated, the GLPK solver package was employed by our students to produce solutions that allow for further analysis.
After formulating this model, our students used the GLPK solver package to produce solutions that allow for further analysis.

%%%%%%%%%%%%%%%%%%%%%%%%%%%%%%%%%%%%%%%%%%%%%%%%%%%%%%%%%%%%%%%%%%%%%%%%%%%%%%%%%%%%%%%%%%%

\subsection{Results}
\label{subseq:6}
For each subset of three or four destination cards, we were able to find an optimal solution that maximizes a player's total points while connecting the destination cities.  Table \ref{T12} shows the percentage of times edges show up in optimal solutions when looking at all three-subsets and four-subsets of destination cards.  

\begin{table}[htbp!] 
\begin{center}
\begin{tabular}{|c|c||c|c|} \hline
 3-subset Edge & Percentage & 4-subset Edge & Percentage    \\ \hline
 \hline
 Duluth--Helena & 6.19\% &  Helena--Seattle & 6.80\%\\ \hline
 Duluth--Sault St. Marie & 5.72\% & Duluth--Helena & 5.36\% \\ \hline
 El Paso--Houston & 5.27\% & Seattle--Vancouver & 4.62\%  \\ \hline
 Helena--Seattle & 5.20\% & Duluth--Toronto & 4.41\% \\ \hline
 El Paso--Los Angeles & 4.52\% & Portland--Seattle & 4.00\% \\ \hline
 Denver--Helena & 3.79\% & Pittsburgh--New York & 3.59\% \\ \hline
 Miami--New Orleans & 3.47\% &  El Paso--Los Angeles & 3.45\%\\ \hline
 Pittsburgh--Toronto & 3.44\% & Nashville--Pittsburgh & 3.44\% \\ \hline
Nashville--Pittsburgh & 3.42\% & Pittsburgh--Toronto & 3.38\%  \\ \hline
 Atlanta--Nashville & 3.42\% & Atlanta--Nashville & 3.32\% \\ \hline
 Houston--New Orleans & 3.23\% & Denver--Santa Fe & 3.14\% \\ \hline
 Portland--Seattle & 3.16\% & El Paso--Houston & 3.08\% \\ \hline
\end{tabular}
\caption{Top edges used in maximizing a player's point total in \emph{Ticket to Ride}.}
\label{T12}
\end{center}
\end{table}

Most of the edges listed in Table \ref{T12} have route length six, each scoring 15 points if claimed during the game.  It is well-known \cite{Slade} that the game over-rewards six-train edges, earning 2.5 points/train, which is $20\%$ more than any other route length.  As seen in Figure \ref{TTR}, these six-train routes form parallel `spines' across the upper and lower game board that can be used to connect east- and west-coast cities, which garner more destination card points.  Further, these spine routes can be branched off of to connect to other destination cities of the subset of cards a player has.  All of these factors, we believe go into bolstering their value in Table \ref{T12}.  We observe that there is one notable non-six-train edge listed among the most used.  The edge between Atlanta and Nashville uses only one train and scores only one point for use, but it is listed in the top ten most used routes in both the three and four-subset scenarios.  We believe this is due to its centrality of routes connecting cities on destination cards.

\resproject{Based off of the inclusion of the Atlanta-to-Nashville edge in many optimal solutions, one interesting project could be to study centrality measures of the \emph{Ticket to Ride} network.  In particular, \emph{betweenness} is a measure of how many shortest paths each edge is a part of.  Routes with high betweenness might be used more often by players.  Studying the betweenness of each edge might help players understand not only their route choices but also route choices that can be used to block other players. See \cite{Newman} for an introduction to betweenness.}
\resproject{\emph{Ticket to Ride} allows players to use the same edge for multiple paths between destination cards.  What if players were required to make edge-disjoint paths?  What would be the optimal routes to take?  Would it be possible to create these for every subset of three destination cards with 45 trains? }
\resproject{\emph{Ticket to Ride} has many variations, including \emph{Ticket to Ride Europe}, \emph{Ticket to Ride San Francisco}, and \emph{Ticket to Ride Rails and Sails} to name a few.  Extend this work to these variants, finding the most popular edges used in forming routes.}
\resproject{In our analysis, we make the assumption that all train allocations should go towards connecting destination card cities.  What if this is not the optimal strategy for a player?  Augment the given IP to make a choice of whether to connect two given destination cities.  Also, it might be useful to use an edge without connecting it to other used edges just to score more points.  Modify the constraints to allow for this.  Use this to answer the question of what is the optimal set of destination cards to have and routes to choose to maximize the number of points scored in a game.}

%%%%%%%%%%%%%%%%%%%%%%%%%%%%%%%%%%%%%%%%%%%%%%%%%%%%%%%%%%%%%%%%%%%%%%%%%%%%%%%

\section{Historical Cryptology}
\label{sec:5}

Historical cryptology is the science and art of making and breaking of codes and ciphers from ancient times up through the 20th century. There are two sides to cryptology: cryptography,  making and implementing codes and ciphers, and cryptanalysis, the breaking of codes and ciphers. To introduce the fundamentals of cryptography, suppose John wants to share a secret message with his friend Liz over an insecure channel (radio, email, etc.).  John and Liz would first decide on a scrambling method (encryption algorithm) and privately share a key that designates exactly how this encryption algorithm will be used. Then John would take his initial message (plaintext) and scramble it with the encryption algorithm using the agreed upon key to create a ciphertext, which looks like a jumbled up mess. John could then send this ciphertext over an insecure  channel to Liz and she could quickly unscramble (decrypt) the message since she knows the key. This unscrambling process, decryption, frequently involves applying the inverse of your encryption function.

To make this process more explicit, we consider a basic example, the shift cipher, which shifts every letter in the plaintext by a fixed amount to obtain the ciphertext.

\begin{example} \label{ex1}
Suppose John wants to send the plaintext ``Knot Theory'' to Liz using the shift cipher. Before this process starts, they agree upon a private key, which is a fixed number $k \in \{0,\ldots, 25\}$ designating how long of a shift to use. For this scenario, let's say $k =11$.  John would first convert his plaintext to numbers ($a=0, \ldots, z=25$) to get the string of (plaintext) numbers 10,13,14,19,19,7,4,14,17,24. Then, John would use addition modulo $26$ as his encryption algorithm via $e_{11}(p) = p + 11$ modulo $26$, for each plaintext letter $p$ to shift every letter 11 units to the right. This produces the string of (ciphertext) numbers 21,24,25,4,4,18,15,25,2,9. Then John converts this string of numbers back to letters to get the ciphertext ``vyzeespzcj'' which he then can send to Liz. Since Liz knows the key $k=11$, she can then convert the ciphertext to numbers and decrypt via the the function $d_{11}(c) = c -11$ modulo $26$, for each ciphertext number $c$.
\end{example}

For the projects discussed in this section, we will be playing the role of the cryptanalyst, where we have intercepted a ciphertext and while we might know the encryption algorithm, we don't know the exact key that was used. For the shift cipher described above, cryptanalysis is a simple task since there are only $25$ non-trivial keys and one could easily check every key by hand. However, with most historical ciphers, the number of keys is astronomically large and one could not even ask a modern computer to check every single key in a reasonable amount of time. Even when it is possible to break such ciphers by more problem-solving oriented methods by hand, it can become an incredibly time consuming and tedious process. This all motivates more modern approaches to cryptanalysis of historical ciphers that leverage building mathematical models that attempt to produce the best possible decipherment. By ``possible decipherment'' we mean applying some choice of a key (not necessarily the correct one) to the ciphertext. For instance, in Example \ref{ex1}, for each $k \in \{0, \ldots, 25\}$, one could apply any $d_{k}(c) = c - k$ modulo $26$ to the given ciphertext to get a possible decipherment. Thus, the number of possible decipherments of a given ciphertext is equal to the number of keys for the cipher in use. Obviously $k=11$ leads to the best possible decipherment in this example, and any other $k$-value modulo $26$ would lead to a decipherment that is $0\%$ accurate. For other types of ciphers, the accuracy of a possible decipherment can exhibit far more variation.

The projects discussed in this section focus on Monoalphabetic Substitution Ciphers (MASCs), though some of our project ideas could be expanded to other historical ciphers, a few of which are discussed as potential research projects. We refer the reader to  \cite{Bauer} and  \cite{DuninSchmeh} for background on historical ciphers, including cryptanalysis techniques. For encryption with an MASC, each letter is replaced with one and only one letter, that is,  if $p$ in the plaintext is substituted with $q$ in the ciphertext at some point, then that substitution consistently occurs between the plaintext and ciphertext, and $q$ in the ciphertext would always decrypt to $p$. The shift ciphers described earlier provide a small sub-class of MASCs. A key for an MASC corresponds with a bijection on the set of alphabet letters or a permutation of the alphabet. As a result, there are $26! \approx 4.03 \times 10^{26}$ possible keys, which makes it unreasonable to brute force. However, the linguistic structure of the underlying language remains in the ciphertext of an MASC. To assist with cryptanlysis, one can then compare frequencies of individual letters (1-grams), adjacent pairs of letters (2-grams), etc. between the ciphertext and the expected behavior of the underlying language, search for common words, like ``and'' and ``the,'' and use other linguistic properties.  

While there are many mathematical modeling approaches to breaking MASCs (see \cite{SabonchiAkay} for a survey on numerous such approaches), we follow and improve upon the methods described in Ravi and Knight \cite{RaviKnight} for creating an IP to model the process of finding a key that maximizes a measure of expectation based on $1$-gram and $2$-gram frequencies. Since the work of Ravi and Knight  was the only example in the literature that the authors could find on using IPs for cryptanalysis of historical ciphers, we believe this is an open area of research, with plenty of opportunities for undergraduate projects. In what follows, we share both the IP model from Ravi and Knight and one of the IP  models developed as part of an undergraduate research project that a subset of the authors supervised.  We use the abbreviations RK to refer to the Ravi and Knight model and UR to refer to the undergraduate research project model. We note that the UR model uses the same decision variables and constraints as the RK model, but the models differ in terms of objective functions.

%%%%%%%%%%%%%%%%%%%%%%%%%%%%%%%%%%%%%%%%%%%%%%%%%%%%%%%%%%%%%%%%%%%%%%%%%%%%%%%%%%%%%

\subsection{Decision Variables and Objective Function}
\label{subseq:7}
\begin{figure}[h!]
\begin{center}
\includegraphics[width = 1.0\textwidth]{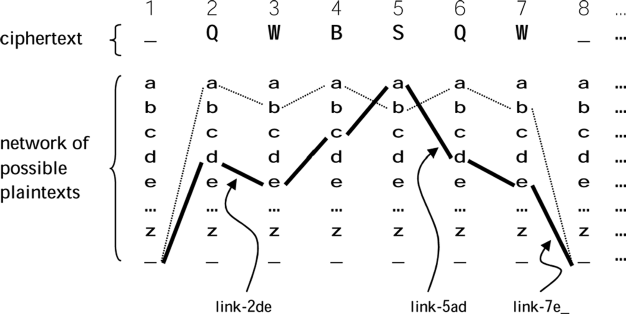}
\caption{A decipherment network highlighting the behavior of linking variables, taken from \cite{RaviKnight}.}
\label{fig:DN}
\end{center}
\end{figure}

 There are two types of binary decision variables introduced and used in the RK model, which are also used in the UR model. Suppose we are given an MASC ciphertext of length $n$ with ordered  ciphertext letters $c_1c_2 \cdots c_n$. The first type of binary variable has a value of $1$ precisely when plaintext letter $p$ is assigned to ciphertext letter $q$.  We will call this binary variable $key_{p,q}$, where $1 \leq p, q \leq 26$  and using a correspondence of $a=1, \ldots, z=26$.   We also set up binary variables $link_{i,j,k}$ to be 1 if ciphertext  letter $c_i$ decrypts to plaintext letter $j$ and $c_{i+1}$ decrypts to plaintext letter $k$, where $1 \leq i \leq n-1$ and $1 \leq j,k \leq 26$. Figure \ref{fig:DN} highlights a visual of the decipherment network from Ravi and Knight where each linking variable corresponds with a path between two consecutive letters in a possible decipherment of the given ciphertext. Here, $``\_"$ designates a space in the text, though for our models, we will assume spaces are not given.  To understand the network, the ordered alphabet, representing potential plaintext letters, has been replicated in columns.  The number of columns $n$ is given by the length of the ciphertext.  Each letter in column $i$, for $2 \leq i \leq n-1$,  is connected via an edge to every letter in columns $i-1$ and $i+1$ although only select edges are shown in Figure \ref{fig:DN}.   Our goal is to find a path through this network of plaintext letters.  The letters in the path will correspond to our possible decipherment of the ciphertext.  In terms of the network in Figure \ref{fig:DN}, the $key_{pq}$ variable serves two roles.  The first is that it assigns ciphertext letter $q$ to plaintext letter $p$.  The second role is similar to the role $y_i^a$ served in the \emph{Ticket to Ride} formulation.  Namely, $key_{pq}$ ``opens" up a path for a $link_{ipr}$ variable to create a path into $p$ and out of $p$ in the $i^{th}$ column of the network.  For the bolded path in Figure \ref{fig:DN}, $key_{dQ}=1$ to allow for $link_{2de}$ and $link_{1\_d}$ to both be 1.  We will describe the constraints that ensure this in Section \ref{subseq:8}.

Ideally, in breaking a cipher, we would like to find the key that maximizes the likelihood that its implementation leads to an accurate decipherment of the ciphertext.  Given an ordered cipheretxt $c_1 c_2\ldots c_n$ and some key $k$ from the set of all possible keys for the given cipher, denoted $\mathcal{K}$, let $p_{1}^{k} p_{2}^{k} \ldots p_{n}^{k}$ represent the decipherment of this ciphertext by applying key $k$.  Then $\{ (p_{1}^{k} p_{2}^{k} \cdots p_{n}^{k})\}_{k \in \mathcal{K}}$ represents the set of all possible decipherments of this ciphertext. Moving forward, we ignore the superscript denoting a particular key to simplify notation.  We seek the possible decipherment $p_1 p_2 \cdots p_n$ that maximizes the probability that this text models the statistical behavior of English, which we denote by $P(p_1 \cdots p_n)$.  Using the laws of conditional probability, we have
$$P(p_1 \cdots p_n) = P(p_n | p_{n-1} \cdots p_1) \cdot P(p_{n-1} | p_{n-2} \cdots p_1) \cdot \ldots \cdot P(p_2 | p_1) \cdot P(p_1 | START),$$ 
where $P(p_1 | START)$ is the probability that the first letter in the plaintext is $p_1$.

In general, finding statistics on $ P(p_n | p_{n-1} \cdots p_1)$ for large $n$ is difficult, so a modeler might choose to approximate this calculation with 
$$P(p_1 \cdots p_n) \approx P(p_n | p_{n-1}) \cdot P(p_{n-1} | p_{n-2}) \cdot \ldots \cdot P(p_2 | p_1) \cdot P(p_1 | START).$$
Several sources exist that give average $n$-gram frequencies for small $n$-values and large sets of text from various languages. See the website Practical Cryptography \cite{Lyons} for the main resource for the UR model. We can use these frequencies to compute the conditional probabilities in the above approximation of $P(p_1 \cdots p_n)$.  Notice that this approximation of $P(p_1 \cdots p_n)$  is multiplicative in nature, but we can transform it to be linear by taking the logarithm of both sides, creating
\[
\begin{array}{lll} \log(P(p_1 \cdots p_n)) & = & \log(P(p_n | p_{n-1})) + \log( P(p_{n-1} | p_{n-2})) + \ldots \\
&&+ \log(P(p_2 | p_1)) + \log(P(p_1 | START)).
\end{array}
\]
However, since the logarithm is negative when its input is between 0 and 1, we multiply each side by $-1$.  The objective function for the RK model \cite{RaviKnight} uses the individual $-\log(P(r|p))$ values to weight the corresponding linking variables:
$$\sum_{i=1}^{n-1}  \sum_{p,r} \bigg(- \log(P(r|p)) \cdot link_{i,p,r}\bigg).$$

For the RK model, we want to minimize the objective function since if $P(r|p)$ has a large probability (better fits behavior of English), then $-\log(P(r|p))$ has a small positive value, while if $P(r|p)$ has a low probability, then $-\log(P(r|p))$ has a large positive value. Ravi and Knight also built and analyzed $1$-gram and $3$-gram models in \cite{RaviKnight}. 

%We now describe the objective function for the UR model. Let $E_{p}$ represent the average relative frequency for letter $p$ in English and let $C_{q}$ represent the relative frequency of letter $q$ in a given ciphertext. Let $C(c_{i+1}|c_{i})$ represent the relative frequency for ciphertext letter $c_{i+1}$ following ciphertext letter $c_{i}$ in the given ciphertext. The objective function (which we want to minimize) for the UR model is given by
%
%$$\sum_{p,q} \bigg(|E_{p} - C_{q}| \cdot key_{p,q}\bigg) + \sum_{i=1}^{n-1} \bigg( \sum_{p,r} |P(r|p) - C(c_{i+1}|c_{i})| \cdot link_{i,p,r} \bigg).$$

We now describe the objective function for the UR model. Where the RK model will choose to turn on linking variables based on the likelihood the letter transition occurs frequently in English, the UR model tries to minimize the difference between the relative frequency use of a transition or a single letter in English versus its frequency in the particular ciphertext.  Let $E_{p}$ represent the average relative frequency for letter $p$ in English and let $C_{q}$ represent the relative frequency of letter $q$ in a given ciphertext. If plaintext letter $p$ encrypts to ciphertext letter $q$, then the expression $|E_p – C_q|$ measures how far away a ciphertext frequency is from the corresponding frequency in English.  Let $C(c_{i+1}|c_{i})$ represent the relative frequency for ciphertext letter $c_{i+1}$ following ciphertext letter $c_{i}$ in the given ciphertext. The expression $|P(r|p) - C(c_{i+1}|c_{i})|$ represents the difference in the ciphertext transition frequency and the corresponding transition frequency occurring in English, if $pr$ encrypts to $c_{i}c_{i+1}$.  The objective function (which we want to minimize) for the UR model is given by
$$\sum_{p,q} \bigg(|E_{p} - C_{q}| \cdot key_{p,q}\bigg) + \sum_{i=1}^{n-1} \bigg( \sum_{p,r} |P(r|p) - C(c_{i+1}|c_{i})| \cdot link_{i,p,r} \bigg).$$
Since we are minimizing, if $key_{p,q} =1$, then its coefficient should be relatively small, which will happen when the relative frequencies of $q$ in the ciphertext behave similarly to the relative frequencies of $p$ in English.  The coefficients on the linking variables serve a similar purpose for transitional probabilities between a pair of letters in the ciphertext and the corresponding pair of letters in the plaintext.

 Recall that there are $26!$ MASC keys, and in general, any interesting cipher is going to have too many keys for a brute force check. Thus, we needed an indirect method for evaluating the quality of a key or an attempted decipherment. That is why neither of the objective functions described above involve weighting every possible MASC key (which is different than the individual variable keys $key_{p,q}$).  However, for the shift cipher, one could build such an objective function. 

\chproblem{Define decision variables and build an appropriate objective function for breaking the shift cipher that evaluates every possible key.}

%%%%%%%%%%%%%%%%%%%%%%%%%%%%%%%%%%%%%%%%%%%%%%%%%%%%%%%%%%%%%%%%

\subsection{Constraints}
\label{subseq:8}

The following set of constraints were first introduced and used in the RK model and then also used in the UR model. 
The first constraint ensures that each plaintext letter maps to exactly one ciphertext letter:
$$\sum_{q=1}^{26} key_{p,q} =1, \hspace{0.1in} 1 \leq p \leq 26.$$
The second constraint ensures that each ciphertext letter gets mapped to by exactly one plaintext letter:
$$\sum_{p=1}^{26} key_{p,q} = 1, \hspace{0.1in} 1 \leq q \leq 26.$$
Collectively, these first two constraints guarantee that our encryption/decryption algorithm is a MASC, though not necessarily the one used to encrypt the given text. 

We also need to form a step-by-step path through the decipherment network, ensuring that the linking variables  are consistent with each other.  The following constraint ensures this:

$$\sum_{p=1}^{26} link_{i,p,r} = \sum_{p=1}^{26} link_{i+1,r,p}, \hspace{0.1in} 1 \leq i \leq n-2, \hspace{0.1in}   1 \leq r \leq 26.$$
Finally, we need to ensure that the chosen linking variables are consistent with the key variables.  To do this, we set up

$$\sum_{r =1}^{26} link_{i,r,p} = key_{p, c_{i+1}}, \hspace{0.1in} 1 \leq i \leq n-1, \hspace{0.1in} 1 \leq p \leq 26,$$

\noindent
where $c_{i+1}$ is the ciphertext letter in the $i+1$ location in the ciphertext. If  for some $r$, we have  $link_{i,r,p} =1$, then this implies that $c_{i+1}$ decrypts to $p$. However, for a MASC, this implies that any other times $q = c_{i+1}$ shows up in the ciphertext, then it must also decrypt to $p$. Hence, we must have $key_{p,c_{i+1}} = 1$ for consistency.

%%%%%%%%%%%%%%%%%%%%%%%%%%%%%%%%%%%%%%%%%%%%%%%%%%%%%%%%%%%%%%%%%%%%%%%%%%%%%%%%%%%%%%%%%%%%%%%

\subsection{Results}
\label{subseq:9}

In Ravi and Knight \cite{RaviKnight},  $1$-gram, $2$-gram, and $3$-gram IP models are evaluated. Table 2 in their paper provides the number of variables, constraints and average run times as a function of cipher lengths $n \in \{8,16,32, 64\}$. Figure 3 in their paper shows the decipherment error for these three models as a function of cipher length $n \in \{2, 4, 8, 16, 32, 64, 128, 256\}$, where $50$ ciphertexts for each such length were tested.   These tests reveal that the $1$-gram model is highly inaccurate. The $3$-gram model is more accurate than the $2$-gram model, though there is eventually a noticeable trade-off in run times. Moving forward, when we refer to the RK model, we specifically mean their $2$-gram model.

\begin{figure}
     \centering
     \begin{subfigure}[b]{0.9\textwidth}
         \centering
         \includegraphics[width=\textwidth]{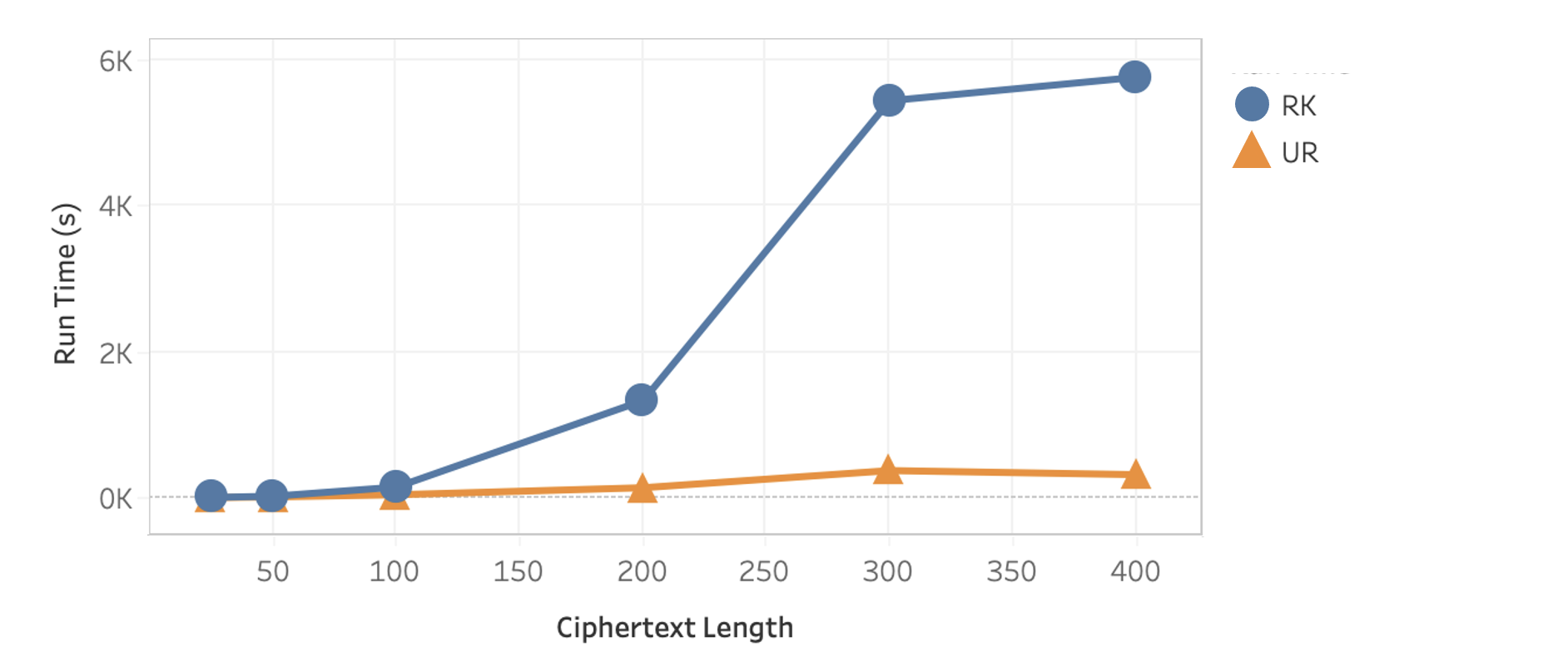}
	\caption{Run time (in seconds).}
	\label{fig:runtime}
     \end{subfigure}
     \begin{subfigure}[b]{0.9\textwidth}
         \centering
         \includegraphics[width=\textwidth]{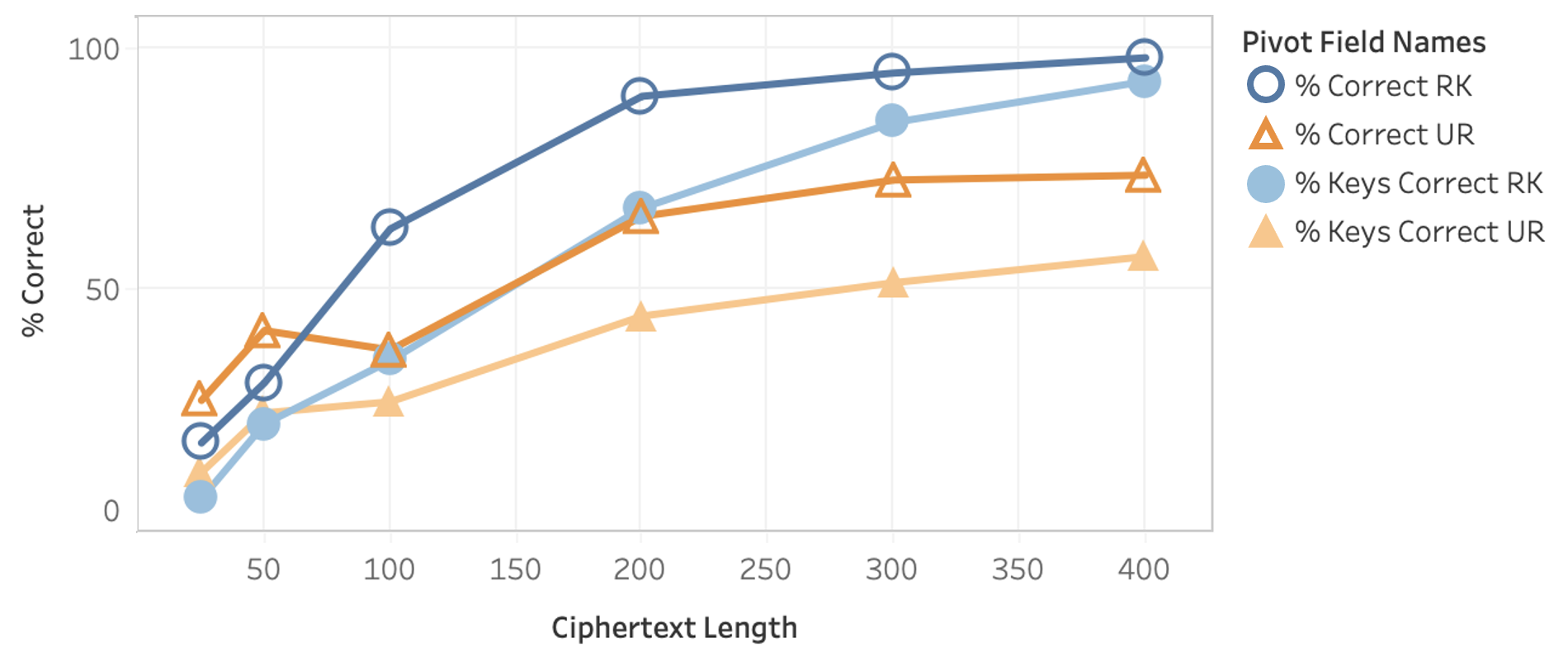}
	\caption{Percentage of correct characters and percentage of correct keys.}
	\label{fig:percentcorrect}
     \end{subfigure}

        \caption{Comparison of the RK and UR models for texts of various lengths.}
        \label{fig:modelcomparison}
\end{figure}

For our UR model, we considered texts with varying lengths, but where the underlying language was always English and no spaces or punctuation were included.  We note that Figure 4 in \cite{RaviKnight} shows the decipherment error for the RK model with spaces and without spaces. Not surprisingly, the model without spaces has a larger decipherment error for nearly every ciphertext length considered. Figure \ref{fig:modelcomparison}  shows some data comparing our UR model with the RK model, which we reproduced. 

Both of these models used a Gurobi solver on a single machine with 64 bit, 2.424 GHz Apple M2 Chip, and 8 GB of memory. For each text length, five texts were used to compare model performance. Figure \ref{fig:runtime} shows the run time in seconds while Figure \ref{fig:percentcorrect} shows both the percentage of letters in the ciphertext that decrypted correctly (as `\% Correct') and the percentage of the $26$ plaintext-ciphertext assignments of an MASC key were correctly identified (as `\% Keys Correct'). 

We  quickly see some tradeoffs between these two models. The UR model has a much faster run time, and this gap in run time only widens as text length increases. The UR model was more accurate on texts of shorter length, while the RK model was more accurate on longer texts. For a more careful analysis, one would want to use several more ciphertexts of each length, and perhaps other ciphertext lengths. 

We also want to emphasize that having a partially correct decipherment that is reasonably accurate can quickly lead to a complete decipherment by hand. We provide a few examples of plaintexts and the  decipherments of the corresponding ciphertexts from the UR model to emphasize this point. We include spaces for the plaintexts and the partial decipherments here for the ease of comparison. Not having spaces makes this process only slightly more challenging with a reasonable partial decipherment. 

\begin{example}
This first text has $395$ characters without punctuation or spaces. The UR model correctly deciphered $91.86\%$ of the ciphertext letters. One can easily recover the rest of the plaintext given this partial decipherment. \\

\textbf{Plaintext:} 
\begin{flushleft}
The Time Traveller for so it will be convenient to speak of him was expounding a recondite matter to us His grey eyes shone and twinkled and his usually pale face was flushed and animated The fire burned brightly and the soft radiance of the incandescent lights in the lilies of silver caught the bubbles that flashed and passed in our glasses Our chairs being his patents embraced and caressed us rather than submitted to be sat upon and there was that luxurious after dinner
\end{flushleft}

\textbf{UR Model decipherment:}

\begin{flushleft}
The time traveller for so it will be zonvenient to sjeap of him was egjoundinc are zondite matter to us his jale crek ekes shone and twinpled and his usuallk jale faze was flushed and animated the fire burnt brichtlk and the soft radianze of the inzandeszent lichts in the lilies of silver zaucht the bubbles that flashed and jassed in our classes our zhairs beinc his jatent sembrazed and zaressed us rather than submitted to be sat ujon and there was that lugurious after dinner
\end{flushleft}
\end{example}

\begin{example}

This second text has  only $90$ characters without punctuation or spaces. The UR model correctly deciphered $57.78\%$ of the ciphertext letters. Even with this level of accuracy, we still see certain words showing up, such as ``the'' and ``she.'' In addition, one can use this partial decipherment to quickly make improvements - ``chullenge'' and ``thut'' can lead to replacing ``u'' with ``a,'' in the decipherment, and so on.  \\

\textbf{Plaintext:} 
\begin{flushleft}
despite the	challenge she persisted	with determination knowing that one	achievement some day awaited her
\end{flushleft}

\textbf{UR Model decipherment:} 
\begin{flushleft}
resmate	the	chullenge she medsaster	jath retedvanutain	onijang	thut ine uchaezevent sive ruq ujuater	hed
\end{flushleft}
\end{example}

Since the RK model and the UR model all use the same decision variables and constraints, it's natural to ask if a different setup could be used to break MASCs.

\resproject{Create an IP model for breaking MASCs that uses a different set of decision variables and constraints. Evaluate and compare with known models. }

As shown in Figure \ref{fig:runtime}, run times quickly increase for the RK model as the ciphertext length increases past $200$ characters. While this isn't a major issue for the UR model, long run times could be an issue for this model for ciphertexts longer than the ones we considered. Thus, one could try to find ``cut-off'' lengths on a ciphertext to determine when the given model does an acceptable job determining most of the key just based on that cut-off portion of the text, while keeping the run time below a given threshold. At the same time, perhaps these models could be tweaked for texts that are sufficiently long so they run faster but still maintain a reasonably high level of accuracy.  This all motivates the following research project.

\resproject{Develop heuristics or modify the RK and UK models to accurately deciphering long ciphertexts while staying within a time threshold.}

%In addition, one could analyze how to modify the UR objective function to improve accuracy across texts of varying length. The UR objective function has two pieces: the sum of the key variables with their coefficients and the sum of the linking variables with their coefficients. A weight $\alpha \in [0,1]$ could be given to the sum of the key variables (with their coefficients) and weight $1 - \alpha$ could be given to the sum of the linking variables (with their coefficients).
%
%\resproject{Modify the UR objective function with weights and test decipherment accuracy for varying text lengths and varying $\alpha$ values. Attempt to explain why weighting one part of the objective function over the other works better for varying text lengths.} 
%
%\christian{Kevin, you already suggested the above research project as one for the students to think about if they want to continue beyond this summer, so I'm not sure if we should include this or not.} \kevin{You're right.  Would you suggest to remove it or replace it with something else?}

There is a vast variety of other historical ciphers that could be attacked using IP models. A natural next step would be to consider homophonic substitution ciphers (HSCs), where each plaintext letter can be replaced with multiple ciphertext letters or symbols, but each ciphertext letter or symbol corresponds with only one plaintext letter. A key for a HSC is given in Table \ref{TableHSC}. There are some similarities to MASCs here: a fixed set of substitution rules are used where each ciphertext symbol corresponds with only one plaintext symbol. However, by using multiple symbols to represent more frequent letters, like e,t,a,i,o, one can mask some of the frequency patterns of the underlying language. For instance, using the key in Table \ref{TableHSC}, we would replace ``e'' in the plaintext with four different values in our ciphertext. 
\begin{table}[htbp!] 
\begin{center}
\begin{tabular}{|c|c|c|c|c|c|c|c|c|c|c|c|c|c|c|c|c|c|c|c|c|c|c|c|c|c|c|} \hline
 Plaintext: & A & B & C & D & E & F & G & H & I & J & K & L & M & N & O & P & Q & R & S & T & U & V & W & X & Y & Z   \\ \hline
 Ciphertext: & 78 & 96 & 26 & 55 & 02 & 37 & 71 & 63 & 04 & 15 & 87 & 29 & 90 & 22 & 81 & 57 & 40 & 16 & 71 & 72 & 67 & 94 & 74 & 10 & 80 & 97 \\ \hline
  Ciphertext: & 52 &  & &  & 29  &  &  & 14 & 81 &  &  &  &  & 11 & 34 &  &  & 35 & 58 & 93 & 21 &  &  & &  &  \\ \hline
  Ciphertext: & 17 &  & &  & 84  &  &  &  & 76 &  &  &  &  & 63 & 21 &  &  & 86 & 29 & 27&  &  &  & &  &  \\ \hline
  Ciphertext: &  &  & &  & 07  &  &  &  &  &  &  &  &  &  &  &  &  & & &  44 &  &  &  & &  &  \\ \hline
\end{tabular}
\caption{A key for a homophonic substitution cipher.}
\label{TableHSC}
\end{center}
\end{table}
\chproblem{Suppose you wanted to create an IP model for breaking the a HSC and planned to use the MASC IP models as a guideline. Further, suppose you know that any plaintext letter maps to at most four distinct ciphertext symbols; the key in Table \ref{TableHSC} meets these qualifications. How many key variables would you have in your model? Write down the appropriate ``key'' constraints for this model, similar to the first two constraints for the MASC models in this section.}

\resproject{Create and evaluate an IP model for breaking HSCs, where the maximum number of ciphertext symbols replacing any plaintext letter is known.}

\section{Conclusion}
Integer programming is a mathematical tool that can be applied to a variety of interesting topics, including those involving decisions. This chapter has introduced three such areas that are accessible to undergraduates without needing much specific experience before starting their research project. We have used board games, such as \emph{The Genius Square} and \emph{Ticket to Ride}, and cryptography as examples of using IP to explore the intricacies of these scenarios. 

The work described in these projects was completed in strong collaboration with undergraduate students working with the authors in various configurations as Summer Mathematics Undergraduate Research Fellows at Furman University. For the project analyzing \emph{The Genius Square} game, students Hayes Brown, Ellis Edinkrah, Joey Maness, and Riley McLachlan worked with faculty mentors E. Bouzarth, J. Harris, and K. Hutson. For the project involving \emph{Ticket to Ride}, students Claire Gillaspy, Justin Hager, Eleanor Liu, Salem Wear, and Max Young collaborated with faculty mentors E. Bouzarth, J. Harris, and K. Hutson. For the cryptography project, students Ellie Johnson and Sam Schaich worked with K. Hutson and C. Millichap. The collaborative efforts of these student-faculty research teams led to productive mathematical explorations of ways integer programming can apply to a variety of interesting subject areas.

\printindex
\end{document}